\documentclass{article}

\usepackage{arxiv}
\usepackage[utf8]{inputenc} 
\usepackage[T1]{fontenc}    
\usepackage[colorlinks]{hyperref}       
\usepackage{url}            
\usepackage{booktabs}       
\usepackage{amsfonts}       
\usepackage{amsmath}       
\usepackage{amssymb}
\usepackage{nicefrac}       
\usepackage{microtype}      
\usepackage{lipsum}
\usepackage{graphicx}
\usepackage{todonotes}
\usepackage{float}
\usepackage[
backend=biber,
style=numeric,
sorting=none,
maxcitenames=99,
maxbibnames=99
]{biblatex}
\usepackage{caption} 
\usepackage{subcaption}

\DeclareMathOperator*{\bigveebar}{\underline{\bigvee}}

\usepackage{xcolor}
\usepackage[linesnumbered,ruled,vlined]{algorithm2e}

\usepackage[affil-it]{authblk}    
\usepackage[toc,page]{appendix}

\SetCommentSty{mycommfont}

\SetKwInput{KwInput}{Input}                
\SetKwInput{KwOutput}{Output}              

\usepackage{tikz}
\definecolor{lime}{HTML}{A6CE39}
\DeclareRobustCommand{\orcidicon}{
	\begin{tikzpicture}
	\draw[lime, fill=lime] (0,0) 
	circle [radius=0.16] 
	node[white] {{\fontfamily{qag}\selectfont \tiny ID}};
	\draw[white, fill=white] (-0.0625,0.095) 
	circle [radius=0.007];
	\end{tikzpicture}
	\hspace{-2mm}
}
\foreach \x in {A, ..., Z}{\expandafter\xdef\csname orcid\x\endcsname{\noexpand\href{https://orcid.org/\csname orcidauthor\x\endcsname}
			{\noexpand\orcidicon}}
}

\renewbibmacro{in:}{%
  \ifentrytype{article}{}{\printtext{}}%
}
\DeclareUnicodeCharacter{0301}{\'{e}}

\title{Logic-Based Discrete-Steepest Descent: A Solution Method for Process Synthesis Generalized Disjunctive Programs}
\shorttitle{Logic-Based Discrete-Steepest Descent: A Solution Method for Process Synthesis GDP}

\author[1,2]{Daniel~Ovalle}
\author[3]{David A.~Li{\~n}\'{a}n\orcidA{}}
\author[4]{Albert Lee\orcidC{}}
\author[2]{Jorge M.~G\'{o}mez}
\author[3]{Luis~Ricardez-Sandoval}
\author[1]{Ignacio E. Grossmann}
\author[1,4,5,6,*]{David E. Bernal Neira\orcidB{}}
\affil[1]{Department of Chemical Engineering, Carnegie Mellon University, Pittsburgh, PA 15213, USA}
\affil[2]{Department of Chemical and Food Engineering, Universidad de los Andes, Bogot\'{a} 1100123, Colombia}
\affil[3]{Department of Chemical Engineering, University of Waterloo, Waterloo, Ontario, Canada N2L 3G1}
\affil[4]{Davidson School of Chemical Engineering, Purdue University, West Lafayette, IN, USA}
\affil[5]{Quantum Artificial Intelligence Laboratory, NASA Ames Research Center, Mountain View, CA, USA}
\affil[6]{Research Institute of Advanced Computer Science, Universities Space Research Association, Mountain View, CA, USA}
\affil[*]{Corresponding author: dbernaln@purdue.edu}

\begin{document}
\maketitle







\begin{abstract}
Optimization of chemical processes is challenging due nonlinearities arising from chemical principles and discrete design decisions.
The optimal synthesis and design of chemical processes can be posed as a Generalized Disjunctive Programming (GDP)  problem.
While reformulating GDP problems as Mixed-Integer Nonlinear Programming (MINLP) problems is common, specialized algorithms for GDP remain scarce.
This study introduces the Logic-Based Discrete-Steepest Descent Algorithm (LD-SDA) as a solution method for GDP problems involving ordered Boolean variables. 
LD-SDA transforms these variables into external integer decisions and uses a two-level decomposition: the upper-level sets external configurations, and the lower-level solves the remaining variables, efficiently exploiting the GDP structure.
 In the case studies presented in this work, including batch processing, reactor superstructures, and distillation columns, LD-SDA consistently outperforms conventional GDP and MINLP solvers, especially as problem size grows. LD-SDA also proves to be superior when solving challenging problems where other solvers encounter difficulties in finding optimal solutions.

\end{abstract}

\keywords{Superstructure Optimization \and Optimal Process Design \and Generalized Disjunctive Programming \and MINLP \and Process Intensification}

\section{Introduction}
The ongoing research in modeling and optimization provides computational strategies to enhance the efficiency of chemical processes across various time scales, e.g., design, control, planning, and scheduling \parencite{grossmann2005enterprise,andres2022integration}. 
Optimization tools assist in developing novel processes and products that align with environmental, safety, and economic standards, thus promoting competitiveness.
Despite advances in the field, deterministically solving optimization problems that involve discrete decisions alongside nonlinearities remains a significant challenge. 

A key application where these challenges arise is the optimal synthesis and design of reactor and separation processes, which requires incorporating discrete decisions to determine the arrangement and sizes of distillation sequences and reactors, as well as the non-ideal relationships needed to model vapor-liquid phase equilibrium. 
The interaction between nonlinear models and discrete decisions in this problem introduces nonconvexities and numerical difficulties (e.g., zero-flows of inactive stages/units), that complicate the direct solution of these problems with traditional optimization solvers \parencite{mencarelli2020review, Ma2021SynthesisDesign}.
Another important area where these challenges become evident is in economic nonlinear model predictive control, where solving optimization problems within the controller sampling time presents a significant obstacle \parencite{CASPARI202050}. 
This difficulty increases when design or scheduling decisions are coupled with control, adding discrete decisions to the formulation \parencite{andres2022integration}.

The computational burden of these problems presents a significant limitation, often preventing the ability to find timely solutions, particularly in online applications or large-scale systems.
This challenge highlights the need for advanced optimization algorithms capable of efficiently exploring the search space of discrete variables while managing nonlinear discrete-continuous interactions to address relevant chemical engineering optimization problems. 
Two major modeling approaches that address these issues by incorporating discrete decisions and activating or deactivating groups of nonlinear constraints are Mixed-Integer Nonlinear Programming (MINLP) and Generalized Disjunctive Programming (GDP).


Typically, optimization problems are posed using MINLP formulations that incorporate both continuous variables, here denoted as $\mathbf{x} = (x_1, \dots, x_{n_x})$, and discrete variables, here denoted as $\mathbf{z} = (z_1, \dots, z_{n_z})$.
The resulting optimization problems involve the minimization of a function $f: \mathbb{R}^{n_x} \times \mathbb{Z}^{n_z} \to \mathbb{R}$ subject to nonlinear inequality constraints $\mathbf{g}: \mathbb{R}^{n_x} \times \mathbb{Z}^{n_z} \to  \mathbb{R}^{l}$.
Variables are usually considered to belong to a closed and bounded set  $\mathbf{x} \in [\underline{\mathbf{x}},\overline{\mathbf{x}}]$ and $\mathbf{z} \in  \{ \underline{\mathbf{z}},\dots,\overline{\mathbf{z}} \}$, respectively.
Throughout this work, we use underbar and overbar notation to denote lower
and upper bounds, respectively.
The mathematical formulation of an MINLP is as follows:
\begin{equation}\label{eqn:minlp.general}
\begin{aligned}
\min_{\mathbf{x},\mathbf{z}} &\ f(\mathbf{x}, \mathbf{z}) \\
\text{s.t.} &\ \mathbf{g}(\mathbf{x}, \mathbf{z}) \leq \mathbf{0}\\
&\ \mathbf{x} \in [\underline{\mathbf{x}},\overline{\mathbf{x}}] \subseteq \mathbb{R}^{n_x}; \mathbf{z} \in \{ \underline{\mathbf{z}},\dots,\overline{\mathbf{z}} \} \subseteq \mathbb{Z}^{n_z}\\
\end{aligned} \tag{MINLP}
\end{equation}

Problem \eqref{eqn:minlp.general} belongs to the NP-hard complexity class \parencite{liberti2019undecidability}.
Despite the challenges associated with solving NP-hard problems, efficient MINLP solution algorithms have been developed, motivated by its various applications \parencite{grossmann2002review}.
These algorithms leverage both discrete and continuous information to find a feasible optimal solution$(\mathbf{x}^*,\mathbf{z}^*)$.
Among the most commonly used approaches to solve \eqref{eqn:minlp.general} deterministically are methods based on decomposition or branch-and-bound (BB) \parencite{lee2011mixed}.
These techniques separately address the two main sources of difficulty in \eqref{eqn:minlp.general}: the discreteness of $\mathbf{z}$ and the nonlinearity of $\mathbf{g}$.

Both decomposition and branch-and-bound rely on bounding the optimal objective value $f(\mathbf{x}^*,\mathbf{z}^*)$.
This involves searching for values ($\hat{f}, \tilde{f} $) such that $\hat{f} \leq f^* \leq \tilde{f}$, and progressively tighten them.
We denote $f^*$ as an optimal solution value, meaning that $f(\mathbf{x}^*,\mathbf{z}^*)=f^*$ if the $(\mathbf{x}^*,\mathbf{z}^*)$ is an optimal solution.
An upper bound on an optimal solution value is found by identifying feasible solutions to the problem  $\{ (\tilde{\mathbf{x}},\tilde{\mathbf{z}}) \mid \mathbf{g}(\tilde{\mathbf{x}},\tilde{\mathbf{z}}) \leq 0, \tilde{\mathbf{x}} \in [\underline{\mathbf{x}},\overline{\mathbf{x}}], \tilde{\mathbf{z}} \in \{ \underline{\mathbf{z}},\dots,\overline{\mathbf{z}} \} \} $, i.e., $f(\mathbf{x}^*,\mathbf{z}^*) \leq f(\tilde{\mathbf{x}},\tilde{\mathbf{z}}) = \tilde{f}$.
Relaxations of problem \eqref{eqn:minlp.general}, which are optimization problems defined over a larger feasible set, have an optimal solution $\hat{f}$ that is guaranteed to underestimate an optimal objective value of the original problem i.e., $\hat{f} \leq f^*$. 

The second modeling approach used in the literature is GDP, which generalizes the problem in \eqref{eqn:minlp.general} by introducing Boolean variables $\mathbf{Y}$ and disjunctions into the formulation
\parencite{trespalacios2014review}. 
In GDP, the Boolean variable $Y_{ik}$ indicates whether a set of constraints
$\mathbf{h}_{ik}(\mathbf{x},\mathbf{z}) \leq \mathbf{0}$ is enforced or not.
We refer to this enforcing condition as a \emph{disjunct} $i$, in \emph{disjunction} $k$. 
Only one disjunct per disjunction is to be selected; hence, we relate disjunctions with an exclusive OR (XOR $\underline{{\vee}}$) operator
, which can be interpreted as an $\texttt{Exactly}(1,\cdot)$ operator when there are more than two disjuncts in a disjunction \parencite{perez2023modeling}.
Boolean variables can be related through a set of logical propositions $\Omega(\mathbf{Y})=\textit{True}$ by associating them through the operators
AND ($\wedge$), OR($\vee$),  XOR ($\underline{{\vee}}$), negation ($\neg$), implication ($\Rightarrow$) and  equivalence ($\Leftrightarrow$).
Furthermore, GDP considers a set of \emph{global} constraints $\mathbf{g}(\mathbf{x}, \mathbf{z}) \leq 0$ existing outside the disjunctions, which are enforced regardless of the values of the Boolean variables.
The mathematical formulation for GDP is as follows:

    

\begin{equation}\label{eqn:gdp.general}
\begin{aligned}
\min_{\mathbf{x},\mathbf{Y},\mathbf{z}} &\ f(\mathbf{x}, \mathbf{z}) \\
\text{s.t.} &\ \mathbf{g}(\mathbf{x}, \mathbf{z}) \leq \mathbf{0}\\
&\ \Omega(\mathbf{Y}) = True \\
&\ \bigveebar_{i\in D_k} \left[
    \begin{gathered}
    Y_{ik} \\
    \mathbf{h}_{ik}(\mathbf{x}, \mathbf{z})\leq \mathbf{0}\\
    \end{gathered}
\right] \quad k \in K\\
&\ \mathbf{x}  \in [\underline{\mathbf{x}},\overline{\mathbf{x}}] \subseteq \mathbb{R}^{n_x}; \mathbf{Y} \in \{False,True\}^{n_y}; \mathbf{z} \in \{ \underline{\mathbf{z}},\dots,\overline{\mathbf{z}} \} \subseteq \mathbb{Z}^{n_z}\\
\end{aligned} \tag{GDP}
\end{equation}
where $n_y = \sum_{k \in K} \lvert D_k \rvert$. 
Moreover, Boolean variables may be associated with empty disjunctions and still appear in the logical propositions $\Omega(\mathbf{Y})$ to model complex logic that does not involve a set of constraints $\mathbf{h}_{ik}(\mathbf{x},\mathbf{z}) \leq \mathbf{0}$.

Throughout this work, we make several assumptions: problem \eqref{eqn:gdp.general} has at least one feasible solution, the search space for continuous, integer, and Boolean variables is bounded, and the objective function remains bounded as well. 
Additionally, the main problem and the subproblems obtained by fixing Boolean configurations satisfy the necessary conditions for standard Nonlinear Programming (NLP) and MINLP algorithms to find a solution. 
Specifically, the functions $f, \mathbf{g}$, and $\mathbf{h}_{ik}$ are assumed to be smooth, with available first and second derivatives when solving NLP subproblems.

Different strategies are available to solve problem \eqref{eqn:gdp.general}.
The traditional approach is to reformulate the problem into an \ref{eqn:minlp.general}, and the two classic reformulations are the big-M reformulation (BM) \parencite{wolsey1999integer, raman1994modelling} and the hull or extended reformulation (HR) \parencite{grossmann2003generalized,balas2018disjunctive}.
However, there are algorithms specifically designed for the GDP framework that exploit the intrinsic logic of the problem.
These tailored algorithms include Logic-Based Outer Approximation (LOA) \parencite{turkay1996logic} and Logic-Based Branch-and-Bound (LBB) \parencite{lee2000new}. 

The GDP framework has recently been used in the optimization of chemical processes.
Some modern applications in process design include co-production plants of ethylene and propylene \parencite{pedrozo2021optimal},
reaction-separation processes \parencite{zhang2018rigorous}, membrane cascades \parencite{ovalle2024optimal}, 
and once-through multistage flash process \parencite{ali2017design}.
Other advances in process synthesis include effective modular process \parencite{chen2019effective}, refrigeration systems \parencite{chen2018synthesis},
and optimization of triple pressure combined cycle power plants \parencite{manassaldi2021optimization}. 
Recently, new solvent-based adhesive products \parencite{cui2018comprehensive} and optimal mixtures \parencite{jonuzaj2018design, jonuzaj2017designing} have been designed using this methodology.
The scheduling of multi-product batch production \parencite{wu2021rolling}, blending operations \parencite{ovalle2023study}, refineries \parencite{li2020multi, su2018integrated},
modeling of waste management in supply chains \parencite{mohammadi2021modeling}, and
multi-period production planning \parencite{lara2018global,rodriguez2017generalized} are some modern applications of the GDP framework in planning and scheduling.
We refer the reader to the review by \textcite{trespalacios2014review} for other developments in GDP applications.

 
In many applications, Boolean variables and disjunctions in GDP formulations represent discrete decisions with intrinsic ordering.
Examples of these ordered decisions include selecting discrete locations (e.g., feed location in a distillation superstructure), determining discrete points in time (such as the starting date of a task in scheduling), or specifying integer values (as seen in the number of units in a design problem, either in parallel or series).
A key characteristic of these problems is that increasing or decreasing the value of these discrete decisions implies an ordered inclusion or exclusion of nonlinear equations from the model.
However, the Boolean variables used in \eqref{eqn:gdp.general} typically do not capture this ordered structure, failing to leverage the potential relationships between successive sets of constraints.

To exploit the ordered structure in optimization problems, a solution strategy was recently introduced in the mixed-integer context to solve MINLP superstructure optimization problems more efficiently.
In this approach, ordered binary variables are reformulated as discrete variables (called \emph{external variables}) to explicitly represent their ordered structure \parencite{linan2020optimalI}.
The reformulated problem is then lifted to an upper-level optimization, where a Discrete-Steepest Descent Algorithm (D-SDA) is applied. 
This algorithm is based on principles from discrete convex analysis, which provides a different theoretical foundation for discrete optimization \parencite{murota1998discrete}.

D-SDA has demonstrated its effectiveness in several MINLP applications, including optimal design of equilibrium \parencite{linan2020optimalII} and rate-based catalytic distillation columns~\parencite{linan2021optimal}.
It proved to be more efficient than state-of-the-art MINLP solvers in both computational time and solution quality. 
Early computational experiments applying D-SDA as a logic-based solver for GDP formulations also showed promising improvements in solution quality and computational efficiency. 
These experiments involved case studies such as reactor network design, rate-based catalytic distillation column design, and simultaneous scheduling and dynamic optimization of network batch processes~\parencite{bernal2022process,linan2024discrete}.



This paper presents the Logic-Based D-SDA (LD-SDA) as a logic-based solution approach specifically designed for GDP problems whose Boolean or integer variables follow an ordered structure.
Our work builds on our previous work in \parencite{bernal2022process} and provides new information on the theoretical properties and details of the computational implementation of the LD-SDA as a GDP solver.
The LD-SDA uses optimality termination criteria derived from discrete convex analysis \parencite{murota1998discrete, favati1990convexity} that allow the algorithm to find local optima not necessarily considered by other MINLP and GDP solution algorithms.
This study also presents new computational experiments that showcase the performance of the LD-SDA compared to the standard MINLP and GDP solution techniques.
The novelties of this work can be summarized as follows: 
\begin{itemize}
\item A generalized version of the external variables reformulation applied to GDP problems is presented, thus extending this reformulation from MINLP to a general class of GDP problems.
\item The proposed framework is more general than previous MINLP approaches, allowing the algorithm to tackle a broader scope of problems.
Through GDP, the subproblems can be either NLP, MINLP, or GDP, instead of the previous framework where only NLP subproblems were supported.
\item An improved algorithm that uses external variable bound verification, fixed external variable feasibility via Feasibility-Based Bounds Tightening (FBBT), verification of already visited configurations, and a reinitialization scheme to improve overall computational time. 
\item  The open-source implementation of the algorithm is generalized for any GDP problem, leading to an automated methodology.
Before executing the LD-SDA, the user only needs to identify the variables to be reformulated as external variables and the constraints that relate them to the problem.
This implementation, formulated in Python, is based on the open-source algebraic modeling language Pyomo \parencite{bynum2021pyomo} and its Pyomo.GDP extension~\parencite{chenpyomo}, and can be found in an openly available GitHub repository\footnote{\url{https://github.com/SECQUOIA/dsda-gdp}}.
\end{itemize}


The remainder of this work is organized as follows. 
\S \ref{sec:background} presents a general background in both GDP solution techniques and in discrete-steepest optimization through discrete convex analysis.
\S \ref{sec:LD-SDA} illustrates the external variable reformulation for Boolean variables.
Furthermore, this section formally describes the LD-SDA and discusses relevant properties and theoretical implications. 
The details of the implementation and the algorithmic enhancements are described in \S \ref{sec:implementation details}.
Numerical experiments were conducted to assess the performance of the LD-SDA across various test cases, including reactor networks, batch process design, and distillation columns with and without catalytic stages. The results of these experiments are detailed in \S \ref{sec:results}.
Finally, the conclusions of the work along with future research directions are stated in \S \ref{sec:conclusions}.


\section{Background}
\label{sec:background}


This section serves two primary objectives. 
Firstly, it provides an introduction to the solution methods employed in Generalized Disjunctive Programming. 
Secondly, it describes the Discrete-Steepest Descent Algorithm along with its underlying theoretical framework, discrete convex analysis.


\subsection{Generalized Disjunctive Programming Reformulations Into MINLP}
\label{sec:minlp_reformulation}

A GDP can be reformulated into an MINLP, enabling the use of specialized codes or solvers that have been developed for MINLP problems \parencite{grossmann2013systematic, balas2018disjunctive}.
In general, the reformulation is done by transforming the logical constraints into algebraic constraints and Boolean variables into binary variables \parencite{trespalacios2014review}.
Moreover, MINLP reformulations handle disjunctions by introducing binary decision variables $\mathbf{y}\in \{0,1\}^{n_y}$, instead of Boolean (\(False\) or \(True\)) variables $\mathbf{Y}$. 
The exclusivity requirement of disjunctions is rewritten as the sum of binary variables adding to one, thus implying that only a single binary variable can be active for every disjunction $k \in K$.

\begin{equation}
\begin{aligned}
\mathbf{Y} \in \{False,True\}^{n_y} &\rightarrow \mathbf{y} \in \{0,1\}^{n_y}\\
\Omega (\mathbf{Y}) = True &\rightarrow \mathbf{A} \mathbf{y} \geq \mathbf{a} \\
\bigveebar_{i\in D_k} [Y_{ik}] \Leftrightarrow \texttt{Exactly}(1,[Y_{ik} \  i \in D_k]) &\rightarrow \sum_{i \in D_k} y_{ik} = 1  \\
\end{aligned} \tag{GDP-MINLP}
\end{equation}

In this section, we describe the two most common approaches to transform a GDP into a MINLP, namely the big-M \eqref{eqn:MINLP.BigM} and the hull reformulations \eqref{eqn:MINLP.Hull}.
Different approaches to reformulating the disjunctions of the GDP into MINLP result in diverse formulations. 
These formulations, in turn, yield distinct implications for specific problem solving~\parencite{bernal2021convex}.

The big-M reformulation uses a large positive constant $M$ into the inequalities to either activate or relax constraints based on the values of the binary variables. 
A constraint being active or redundant is dependent on the values taken by their corresponding binary variables.
When a binary variable is $True$, the corresponding vector of constraints \(\mathbf{h}_{ik}(\mathbf{x}, \mathbf{z})\leq \mathbf{0}\) is enforced. 
Otherwise, the right-hand side is relaxed by the large value $M$ such that the constraint is satisfied irrespective of the values of $\mathbf{x}$ and $ \mathbf{z}$, effectively making the constraint nonbinding.
This behavior can be expressed as  \(\mathbf{h}_{ik}(\mathbf{x}, \mathbf{z})\leq M(1-y_{ik})\) where \(y_{ik}\) is a binary variable replacing $Y_{ik}$. 
The resulting GDP-transformed MINLP using \eqref{eqn:MINLP.BigM} is given by:
\begin{equation}\label{eqn:MINLP.BigM}
\begin{aligned}
\min_{\mathbf{x},\mathbf{y},\mathbf{z}} &\ f(\mathbf{x}, \mathbf{z}) \\
\text{s.t.} &\ \mathbf{g}(\mathbf{x}, \mathbf{z}) \leq \mathbf{0}\\
&\ \mathbf{A} \mathbf{y} \geq \mathbf{a} \\
&\ \sum_{i \in D_k} y_{ik} = 1 \quad k \in K \\
&\ \mathbf{h}_{ik}(\mathbf{x}, \mathbf{z}) \leq M_{ik}(1 - y_{ik}) \quad i \in D_k, k \in K\\
&\ \mathbf{x}  \in [\underline{\mathbf{x}},\overline{\mathbf{x}}] \subseteq \mathbb{R}^{n_x}; \mathbf{y} \in \{0,1\}^{n_y}; \mathbf{z} \in \{ \underline{\mathbf{z}},\dots,\overline{\mathbf{z}} \} \subseteq \mathbb{Z}^{n_z}\\
\end{aligned} \tag{BM}
\end{equation}

The hull reformulation (HR) offers an alternative by handling disjunctive inequalities using binary variables, but it takes a different approach by disaggregating both continuous and discrete variables.
For each disjunct in the GDP, a copy of each variable \({v}_{ik}  \in [\underline{{x}},\overline{{x}}]\) or \({w}_{ik} \in \{ \underline{{z}},\dots,\overline{{z}} \} \) 
is created for each element $i$ in the disjunction $D_k$. 
When the corresponding binary variable is set to zero, these new variables become zero as well, effectively deactivating their associated constraints. 
Conversely, only the copies of variables linked to binary variables equal to one are involved in enforcing the constraints. Additionally, the constraints in each disjunct are governed by the binary variables through a perspective reformulation applied to the disaggregated variables. 
Specifically, each disjunct that activates a set of constraints \(\mathbf{h}_{ik}(\mathbf{x}, \mathbf{z})\leq 0\) is reformulated as  \(y_{ik} \mathbf{h}_{ik}(\mathbf{v}_{ik}/y_{ik}, \mathbf{w}_{ik}/y_{ik})\leq 0\), where \(y_{ik}\) is a binary variable that replaces \(Y_{ik}\).

The difficulty in applying the HR to a GDP is that the perspective function \(y_{ik} \mathbf{h}_{ik}(\mathbf{v}_{ik}/y_{ik}, \mathbf{w}_{ik}/y_{ik})\) is numerically undefined when \(y_{ik}=0\) if the constraints in the disjuncts are nonlinear. 
Therefore, the method can potentially fail to find a solution to the GDP problem. 
This issue can be overcome by approximating the perspective function with an inequality as demonstrated in \parencite{furman2020computationally}. 
The resulting GDP-transformed MINLP using \eqref{eqn:MINLP.Hull} yields:
\begin{equation}\label{eqn:MINLP.Hull}
\begin{aligned}
\min_{\mathbf{v},\mathbf{w},\mathbf{x},\mathbf{y},\mathbf{z}} &\ f(\mathbf{x}, \mathbf{z}) \\
\text{s.t.} &\ \mathbf{g}(\mathbf{x}, \mathbf{z}) \leq \mathbf{0}\\
&\ \mathbf{A} \mathbf{y} \geq \mathbf{a} \\
&\ \sum_{i \in D_k} y_{ik} = 1 \quad k \in K \\
&\ \mathbf{x} = \sum_{i \in D_k} \mathbf{v}_{ik} \quad k \in K \\
&\ \mathbf{z} = \sum_{i \in D_k} \mathbf{w}_{ik} \quad k \in K \\
&\ y_{ik}  \mathbf{h}_{ik}(\mathbf{v}_{ik}/y_{ik}, \mathbf{w}_{ik}/y_{ik})\leq 0 \quad i \in D_k, k \in K \\
&\ y_{ik} \underline{\mathbf{x}} \leq \mathbf{v}_{ik} \leq y_{ik} \overline{\mathbf{x}}  \quad i \in D_k, k \in K \\
&\ y_{ik} \underline{\mathbf{z}} \leq \mathbf{w}_{ik} \leq y_{ik} \overline{\mathbf{z}}  \quad i \in D_k, k \in K \\
&\ \mathbf{v}_{ik}  \in [\underline{\mathbf{x}},\overline{\mathbf{x}}] \subseteq \mathbb{R}^{n_x}; \mathbf{w}_{ik} \in \{ \underline{\mathbf{z}},\dots,\overline{\mathbf{z}} \} \subseteq \mathbb{Z}^{n_z}\\
&\ \mathbf{x}  \in [\underline{\mathbf{x}},\overline{\mathbf{x}}] \subseteq \mathbb{R}^{n_x};  \mathbf{y} \in \{0,1\}^{n_y}; \mathbf{z} \in \{ \underline{\mathbf{z}},\dots,\overline{\mathbf{z}} \} \subseteq \mathbb{Z}^{n_z} \\
\end{aligned} \tag{HR}
\end{equation}

The hull reformulation introduces a larger number of constraints compared to the big-M method. 
However, it yields a tighter relaxation in the continuous space, refining the representation of the original GDP problem. 
This can potentially reduce the number of iterations required for MINLP solvers to reach an optimal solution.
Depending on the solver and the problem, the trade-off between these two reformulations might result in one of them yielding problems that are solved more efficiently~\parencite{bernal2021convex}.

Although MINLP reformulations are the standard approach for solving GDP problems, they often introduce numerous algebraic constraints. 
Some of these constraints may be irrelevant to a specific solution and can cause numerical instabilities when their corresponding variables are zero.
This can increase the complexity of the problem and solution time, highlighting the potential of alternative GDP solution techniques that avoid transforming the problem into an MINLP.

\subsection{Generalized Disjunctive Programming Logic-Based Solution Algorithms} 

Instead of reformulating the GDP into an MINLP, and solving the problem using MINLP solvers, some methods developed in the literature aim to directly exploit the logical constraints inside the GDP. 
Attempts to tackle the logical propositions for solving the GDP problem are known as logic-based methods.
Logic-Based solution methods are generalizations of MINLP algorithms that apply similar strategies to process Boolean variables to those used for integer variables in MINLP solvers.
This category of algorithms includes techniques such as Logic-Based Outer-Approximation and Logic-Based Branch-and-Bound \parencite{trespalacios2014review}. 

In GDP algorithms, the (potentially mixed-integer) nonlinear programming subproblems generated upon setting specific discrete combinations, which now encompass logical variables, are confined to only those constraints relevant to the logical variables set to $True$ in each respective combination.
In logic-based algorithms, the generated NLP subproblems, that could potentially be mixed-integer as well, arise from fixing specific Boolean configurations. 
These configurations constrain the (MI)NLP subproblems to only relevant constraints corresponding to logical variables set to $True$ in each setting.
Specifically, when considering a given assignment for the logical variables denoted by \(\hat{\mathbf{Y}}\), the resulting subproblem is defined as:

\begin{equation}
\label{eqn:Sub}
\begin{aligned}
 \min_{\mathbf{x},\mathbf{z}} &\ f(\mathbf{x}, \mathbf{z}) \\
\text{s.t.} &\ \mathbf{g}(\mathbf{x}, \mathbf{z}) \leq \mathbf{0}\\
&\ \mathbf{h}_{ik}(\mathbf{x}, \mathbf{z}) \leq \mathbf{0}\ ~\text{if}~ \hat{{Y}}_{ik} = True \quad i \in D_k, k\in K\\
&\ \mathbf{x}  \in [\underline{\mathbf{x}},\overline{\mathbf{x}}] \subseteq \mathbb{R}^{n_x},\
\ \mathbf{z} \in \{ \underline{\mathbf{z}},\dots,\overline{\mathbf{z}} \} \subseteq \mathbb{Z}^{n_z}\\
\end{aligned} \tag{Sub}
\end{equation}

The formulation of Problem \eqref{eqn:Sub} represents the optimization problem under the constraints governed by the chosen logical assignment \(\hat{\mathbf{Y}}\).
In the most general case, after fixing all Boolean variables, the Problem \eqref{eqn:Sub} is a MINLP. 
Still, in most applications, where there are no discrete decisions besides the ones represented in the Boolean space \( n_{z}=0 \), Problem \eqref{eqn:Sub} becomes an NLP. 
This problem avoids evaluating numerically challenging nonlinear equations whenever their corresponding logical variables are irrelevant (i.e., ``zero-flow'' issues) \parencite{mencarelli2020review}.
The feasibility of Boolean variables in the original equation \eqref{eqn:gdp.general} depends on logical constraints \( \Omega({\hat{\mathbf{Y}}}) = True \).
By evaluating these logical constraints, infeasible Boolean variable assignments can be eliminated without needing to solve their associated subproblems. 

In general, logic-based methods can be conceptualized as decomposition algorithms. At the upper-level problem, these methods focus on identifying an optimal logical combination 
\(\hat{\mathbf{Y}}\). 
This combination ensures that the subproblems \eqref{eqn:Sub}, when solved, converge to an optimal solution of Eq. \eqref{eqn:gdp.general}.
Overall, given a Boolean configuration, the subproblem \eqref{eqn:Sub} is a reduced problem that only considers relevant constraints, and is therefore numerically more stable, and yields faster evaluations than a monolithic MINLP.
Consequently, unlike mixed-integer methods, logic-based approaches can offer advantages, given they exploit the structure of the logical constraints.

{\color{teal}A prevalent logic-based method is the Logic-Based Outer Approximation algorithm, that uses linear relaxations of nonlinear functions to approximate the feasible region of the original problem. 
This approach simplifies complex optimization problems by making them easier to solve while providing bounds on an optimal solution.}
By utilizing linear approximations at iterations \(l=1, \dots, L\) and iterations \(L_{ik} = \{l \mid Y_{ik} = True \text{ for iteration } l \}\), LOA leads to the formulation of a linearized GDP, where an optimal solution provides the integer combinations necessary for problem resolution.
The upper-level problem \eqref{eqn:MainLOA} in the LOA method is as follows:


\begin{equation}
\label{eqn:MainLOA}
\begin{aligned}
 \min_{\mathbf{x},\mathbf{z}, \alpha} &\ \alpha \\
\text{s.t.} &\ \alpha \geq \bar{f}(\mathbf{x},\mathbf{z};\mathbf{x}^l,\mathbf{z}^l) \quad \forall \ l = 1,\dots,L
\\
&\ \bar{\mathbf{g}}(\mathbf{x},\mathbf{z};\mathbf{x}^l,\mathbf{z}^l) \leq 0 \quad \forall \ l = 1,\dots,L\\
&\ \bigveebar_{i\in D_k} \left[
    \begin{gathered}
    Y_{ik} \\
    \bar{\mathbf{h}}_{ik}(\mathbf{x},\mathbf{z};\mathbf{x}^l,\mathbf{z}^l) \leq \mathbf{0} \quad l \in L_{ik} \\
    \end{gathered}
\right] \quad k \in K  \\
&\ \Omega(\mathbf{Y})= True \\
&\ \mathbf{x}  \in [\underline{\mathbf{x}},\overline{\mathbf{x}}] \subseteq \mathbb{R}^{n_x},\
\ \mathbf{z} \in \{ \underline{\mathbf{z}},\dots,\overline{\mathbf{z}} \} \subseteq \mathbb{Z}^{n_z}, \alpha \in \mathbb{R}_+\\
\end{aligned} \tag{Main l-GDP}
\end{equation}

where \(\bar{f}(\mathbf{x},\mathbf{z};\mathbf{x}^l,\mathbf{z}^l)\) is the linear relaxation of function \(f(\mathbf{x},\mathbf{z})\) at point \(\mathbf{x}^l,\mathbf{z}^l )\).
A similar definition is given for the linear relaxations of the global constraints \(\bar{\mathbf{g}}(\mathbf{x},\mathbf{z};\mathbf{x}^l,\mathbf{z}^l)\), and of the constraints inside of the disjunctions \(\bar{\mathbf{h}}_{ik}(\mathbf{x},\mathbf{z};\mathbf{x}^l,\mathbf{z}^l)\).
Inspired by the outer-approximation algorithm for MINLP \parencite{duran1986outer}, these linear relaxations can be built using a first-order Taylor expansion around point \(\{\mathbf{x}^l,\mathbf{z}^l \}\), i.e.,  \(\bar{f}(\mathbf{x},\mathbf{z};\mathbf{x}^l,\mathbf{z}^l) = f(\mathbf{x}^l, \mathbf{z}^l) + \nabla_{x} f(\mathbf{x}^l, \mathbf{z}^l)^\top (\mathbf{x} - \mathbf{x}^l) + \nabla_{z} f(\mathbf{x}^l, \mathbf{z}^l)^\top (\mathbf{z} - \mathbf{z}^l)\).
It is important to note that linear approximations are guaranteed to be relaxations only when the functions $f, \mathbf{g}$, and $\mathbf{h}_{ik}$ are convex. 
For convex nonlinear functions, these linear approximations provide valid bounds on the optimal solution \parencite{kronqvist2019review}.


Problem \eqref{eqn:MainLOA} is usually reformulated into a Mixed-Integer Linear Programming (MILP) problem using the reformulations outlined in \S \ref{sec:minlp_reformulation}.
Upon solving the main MILP problems, the logical combination \(\hat{\mathbf{Y}}\) is determined, defining the subsequent Problem \eqref{eqn:Sub} with the resulting logical combination.
Expansion points for additional constraints are then provided to solve the subproblem \eqref{eqn:Sub} within the context of \eqref{eqn:MainLOA}. While \eqref{eqn:MainLOA} yields a rigorous lower bound, the \eqref{eqn:Sub} subproblem provides feasible solutions, thus establishing feasible upper bounds. 
Each iteration refines the linear approximation of \eqref{eqn:MainLOA}, progressively tightening the constraints and guiding the current solution towards an optimal solution.

Gradient-Based linearizations provide a valid relaxation for convex nonlinear constraints, but do not guarantee an outer-approximation for nonconvex ones.
This limitation jeopardizes the convergence guarantees to globally optimal solutions of LOA for nonconvex GDP problems.
To address this problem, if the linearization of the functions defining the constraints is ensured to be a relaxation of the nonlinear constraints, LOA can converge to global solutions in nonconvex GDP problems.
These relaxations remain to be linear constraints, often constructed using techniques such as multivariate McCormick envelopes \parencite{tsoukalas2014multivariate}.
This generalization is known as Global Logic-Based Outer Approximation (GLOA).


Another important logic-based solution method is the Logic-Based Branch-and-Bound algorithm that systematically addresses GDP by traversing Boolean variable values within a search tree. 
Each node in this tree corresponds to a partial assignment of these variables.
LBB solves optimization problems by splitting them into smaller subproblems with fixed logic variables and eliminating subproblems that violate the constraints through a branch-and-bound technique.  The core principle of LBB is to branch based on the disjunction, enabling it to neglect the constraints in inactive disjunctions.
Furthermore, LBB accelerates the search for an optimal solution by focusing solely on logical propositions that are satisfied. 
Initially, all disjunctions are unbranched, and we define this set of unbranched disjunctions as $KN$.
The LBB starts with the relaxation of the GDP model \eqref{eqn:node GDP} in which all nonlinear constraints from the disjunctions are ignored.
For every node \(l\), the set of branched disjunctions $KB$ can be defined as \(KB^l = K \setminus KN^l \).

\begin{equation}
\label{eqn:node GDP}
\begin{aligned}
    \min_{\mathbf{x},\mathbf{Y},\mathbf{z}} &\ f^l(\mathbf{x}, \mathbf{z}) \\
\text{s.t.} &\ \mathbf{g}(\mathbf{x}, \mathbf{z}) \leq \mathbf{0}\\
&\ \Omega(\mathbf{Y}^l) = True \\
&\ \mathbf{h}_{ik}(\mathbf{x}, \mathbf{z}) \leq \mathbf{0} \quad \text{if  } \hat{Y}^l_{ik} = True, \quad k \in KB^l \\
&\ \bigveebar_{i\in D_k} \left[
    \begin{gathered}
    Y_{ik} \\
    \mathit{\Psi}_{ik}(\hat{\mathbf{Y}}^l) = True \\
    \end{gathered}
\right] \quad k \in KN^l\\
&\ \mathbf{x}  \in [\underline{\mathbf{x}},\overline{\mathbf{x}}] \subseteq \mathbb{R}^{n_x}; \mathbf{Y} \in \{False,True\}^{n_y}; \mathbf{z} \in \{ \underline{\mathbf{z}},\dots,\overline{\mathbf{z}} \} \subseteq \mathbb{Z}^{n_z}\\
\end{aligned}\tag{node-GDP}
\end{equation}

where $\Psi$ denotes the set-valued function of constraints relevant for the unbranched nodes \(KN^l\).

At each iteration, the algorithm selects the node with the minimum objective solution from the queue. 
The objective value of each evaluated node in the queue serves as a lower bound for subsequent nodes. 
Eventually, the minimum objective value among all nodes in the queue establishes a global lower bound on the GDP. 
Branching out all the disjunctions, the algorithm terminates if the upper bound to the solution, determined by the best-found feasible solution, matches the global lower bound.

As mentioned above, logic-based methods leverage logical constraints within the GDP by activating or deactivating algebraic constraints within logical disjunctions during problem solving. 
In the branching process, infeasible nodes that violate logical propositions may be found. These nodes are pruned if they do not satisfy the relevant logical constraints \(\Psi(\mathbf{Y}) = True\).

While logic-based methods offer several advantages, they also have some limitations. 
For nonconvex GDP problems, LOA may struggle to identify the global optimum, as solutions to the NLP subproblems might not correspond to the global optimum. 
Similarly, the LBB method can be resource-intensive, requiring substantial computational time and resources, particularly for large and complex problems \parencite{chenpyomo}.
More specifically, as the problem size increases, the number of subproblems tends to grow exponentially. Thus, there is an ongoing need for more efficient logic-based algorithms that can effectively leverage the logical structure of GDP problems.

The methods described in this section require access to the original GDP problem.
Such an interface has been provided by a few software packages, including Pyomo.GDP \parencite{chenpyomo}.
The LOA, GLOA, and LBB algorithms are evaluated in this work through their implementation in the GDP solver in Pyomo, GDPOpt \parencite{chenpyomo}.


\subsection{Discrete Convex Analysis and the Discrete-Steepest Descent Algorithm} \label{sec:discrete_convex}

Unlike traditional MINLP and GDP solution strategies, which rely on conventional convexity theory treating discrete functions as inherently nonconvex, the Discrete-Steepest Descent Algorithm incorporates an optimality condition based on discrete convex analysis. 
This framework provides an alternative theoretical foundation for discrete optimization, defining convexity structures for discrete functions \parencite{murota1998discrete}.
In this context, a solution to an Integer Programming (IP) problem is considered locally optimal when the discrete variables yield the lowest objective value within a defined \emph{neighborhood}. 
Specifically, this means that the point 
$\mathbf{z}$ must have an objective value lower or equal than all its neighboring points. 
Formally, the neighborhood $N_k$ of a point $\mathbf{z}$ is defined as all the integer points $\mathbf{\alpha}$ (called \emph{neighbors}) within a $k$-ball of radious one centered around $\mathbf{z}$:
\begin{equation}
    \begin{aligned}
        N_k(\mathbf{{z}}) = \{\mathbf{\alpha} \in \mathbb{Z}^{n_{z}}: \; 
\lVert \mathbf{\alpha} - \mathbf{{z}} \rVert_k \leq 1\}
    \end{aligned}
    \label{eq:neighborhood}
\end{equation}

Once the neighborhood of a point $\mathbf{z}$ is identified, the set of directions $\Delta_k(\mathbf{{z}})$ to each of its neighboring points can be computed through vector subtraction, as shown by the following equation,
\begin{equation}
    \begin{aligned}
        \Delta_k(\mathbf{z}) = \{\mathbf{d} : \; \mathbf{\alpha} - \mathbf{{z}} = \mathbf{d}, \forall \ \mathbf{\alpha} \in N_k(\mathbf{{z}})\}
    \end{aligned}
    \label{eq:directions}
\end{equation}
which measures of how far apart the neighbors are in the lattice.

\begin{figure}[ht]
     \centering
     \begin{subfigure}[b]{0.4\textwidth}
         \centering
         \includegraphics[width=\textwidth]{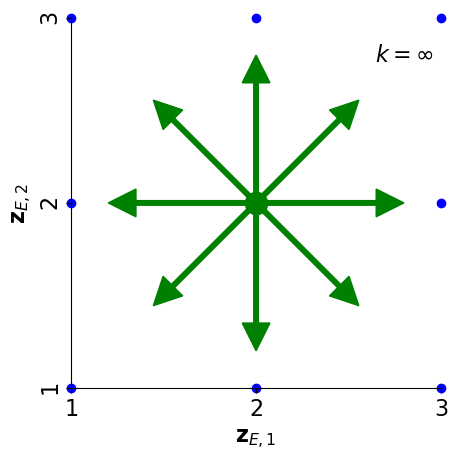}
         \caption{$\infty$-neighborhood ($N_\infty$)}
         \label{fig: KIsInfinity}
     \end{subfigure}
     \hfill
     \begin{subfigure}[b]{0.4\textwidth}
         \centering
         \includegraphics[width=\textwidth]{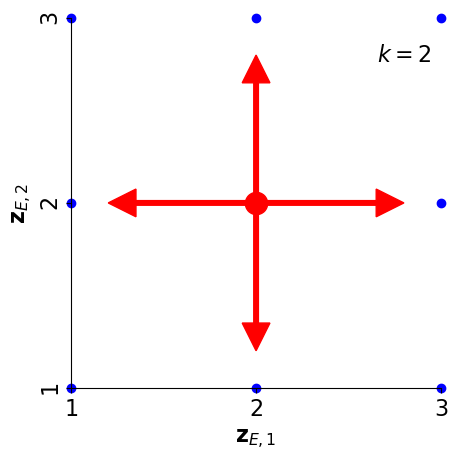}
         \caption{$2$-neighborhood ($N_2$)}
         \label{fig: KIs2}
     \end{subfigure}
        \caption{Visualization of the two neighborhoods $N_\infty$ and $N_2$ on a two-variable discrete lattice, centered at the point $\mathbf{z_E}=(2,2)$. The $\infty$-neighborhood allows movement to all points within unitary $\ell_\infty$ distance, offering a more flexible search space, while the $2$-neighborhood restricts movement to orthogonal directions, providing a more constrained search. This illustrates how the choice of neighborhood affects the directions explored during optimization.}
        \label{fig:neighborhoods}
\end{figure}

In this work, we consider $k \in \{2, \infty\}$ and Figure~\ref{fig:neighborhoods} illustrates both neighborhoods for the case of two dimensions. Therefore, the choice of neighborhood directly affects the local optimum obtained. Importantly, under certain conditions, local optimality within specific neighborhoods can imply global optimality. For example, global optimality is guaranteed for unconstrained integer problems with a separable convex objective function when the positive and negative coordinates of the axes are used as neighbors.

Within the discrete convex analysis framework, an important concept is the notion of \emph{integrally convex} objective functions, as introduced by \textcite{favati1990convexity}.
A function is classified as integrally convex when its local convex extension is convex.
This extension is constructed by linearly approximating the original function within unit hypercubes of its domain, allowing for a more flexible approach to defining convexity in discrete spaces (see \textcite{murota1998discrete}). 
 Integrally convex functions are particularly relevant because they encompass many of the discrete convex functions commonly studied in the literature, including separable convex functions \parencite{murota2023recent}.
 
An MINLP problem may exhibit integral convexity even if it is nonconvex by traditional MINLP standards, where a problem is considered convex if its continuous relaxation is convex \parencite{linan2024discrete}.
This distinction is important for understanding how integral convexity can be leveraged in optimization problems that might otherwise be categorized as nonconvex.
In discrete convex analysis, an optimal solution over the  $\infty$-neighborhood (denoted $N_\infty$) is referred to as \emph{integrally local} (\textit{i-local}),
meaning it is globally optimal for an integrally convex objective function.
Similarly, an optimal solution within $2$-neighborhood (denoted $N_2$), which is sometimes referred to as the separable neighborhood, is known as \emph{separable local} (\textit{s-local}) as it represents a global optimum for a separable convex objective function \parencite{linan2020optimalI,linan2021optimal}.
This distinction between $i$-local and $s$-local optimality is important for understanding how different neighborhoods affect the global optimality guarantees for discrete optimization problems under the lens of discrete convex analysis.
    
The theory of discrete convex analysis was first extended to MINLP in \parencite{linan2020optimalI}. 
Here, the authors introduced a decomposition approach designed for problems with ordered binary variables, which were reformulated using the \emph{external variable} method.
In this approach, external variables were decoupled from the main problem and addressed in an upper-level problem.
Then, the D-SDA was utilized to optimize these external variables, with binary variables fixed accordingly. Furthermore, the objective function values came from the solution of NLP optimization subproblems. 

A significant advantage of using D-SDA is that binary variables reformulated as external variables are evaluated only at discrete points. 
This avoids the issue of evaluating fractional solutions (e.g., a binary variable evaluated at 0.5) given that evaluating points at discrete points suffices to assess discrete optimality requirements.
As a result, D-SDA avoids the potential nonconvexities introduced by the continuous relaxation of MINLP superstructures, e.g., the multi-modal behavior found when optimizing the number of stages in a catalytic distillation column or the number of reactors in series \parencite{linan2020optimalI, linan2020optimalII}.
Similarly, this algorithm was successfully applied to highly nonlinear MINLP problems such as the optimal design of rate-based and dynamic distillation systems \parencite{linan2021optimal,linan2022optimal}. 

Guaranteeing global optimality for MINLPs or GDPs problems remains challenging. 
Unlike previously studied IP problems, where global optimality can sometimes be ensured, MINLP problems involve nonlinear objective functions for each discrete point. This complexity makes it difficult to guarantee global optimality \parencite{murota2023recent,frank2022decreasing}.
The D-SDA seeks to find the best possible solution by choosing appropriate neighborhoods, with the 
$\infty$-neighborhood often used for local optimality in MINLP problems. 
This neighborhood includes all discrete points within an infinity norm of the evaluated point, offering more comprehensive coverage than just positive and negative coordinates, as illustrated in Figure~\ref{fig:neighborhoods}.
Additionally, when applying the D-SDA with the $\infty$-neighborhood to a binary optimization problem (without reformulation), a complete enumeration over discrete variables is required. While not computationally efficient, this method offers a ``brute-force'' alternative for addressing small-scale discrete optimization problems.

Building on these advancements, this paper extends the D-SDA methodology to address more general GDP problems. This extension involves exploring the search space of reformulated Boolean variables directly, eliminating the need for a \eqref{eqn:MINLP.BigM} or \eqref{eqn:MINLP.Hull} reformulation step and avoiding the linearization required in the LOA method.

\section{The Logic-Based Discrete-Steepest Descent Algorithm as a Generalized Disjunctive Programming Algorithm}
\label{sec:LD-SDA}

This section presents Logic-Based D-SDA as a GDP algorithm.
It begins with an explanation of the reformulation process for Boolean variables into external variables, outlining the requirements necessary for reformulation. 
For this, we provide a comprehensive example for demonstration.
Second, the basis of the LD-SDA as a decomposition algorithm that utilizes the structure of the external variables is elucidated. 
The following subsection describes the different algorithms that compose the Logic-Based D-SDA in the context of solving a GDP problem.
The properties of LD-SDA are explained in the final subsection.

\subsection{GDP Reformulations Using External Variables}
\label{sec:reformulation}

Consider GDP problems where a subset of the Boolean variables in $\mathbf{Y}$ can be reformulated into a collection of integer variables referred to as \emph{external variables}.
Thus, $\mathbf{Y}$ in \eqref{eqn:gdp.general} is defined as $\mathbf{Y}=(\mathbf{Y_R},\mathbf{Y_N})$, where $\mathbf{Y_R}=(\mathbf{Y_R}_1,\mathbf{Y_R}_2,...,\mathbf{Y_R}_{n_{R}})$ contains those vectors of independent Boolean variables that can be reformulated using external variables. This means that each vector $\mathbf{Y_R}_j$ will be reformulated with one external variable, and this reformulation is applied for every $j$ in $\{1,2,...,n_R\}$.
Finally, as a requisite to apply the reformulation, each vector $\mathbf{Y_R}_j
 \forall \ j \in \{1,2,...,n_R\}
$ must satisfy the following conditions:   
\begin{itemize}
    \item \textbf{Requirement 1}: Every Boolean variable in $\mathbf{Y_R}_j$ must be defined over a finite well-ordered set $S_j$~\parencite[p.~38]{ciesielski1997set}. This set may be different for each vector of variables; thus, it is indexed with $j$. In addition, variables defined over $S_j$ must represent ordered decisions such as finding discrete locations, selecting discrete points in time, counting the number of times a task is performed, etc.
    Notably, these independent Boolean variables can have indices other than the ordered set.
    Also, not every Boolean variable defined over $S_j$ is necessarily required to be in $\mathbf{Y_R}_j$. For instance, the Boolean variables that determine the feed stage in a distillation column are defined over the set of trays, but some trays may be excluded from $\mathbf{Y_R}_j$ if needed. 
    \item \textbf{Requirement 2}: Boolean variables $\mathbf{Y_R}_j$ are subject to a partitioning constraint $\texttt{Exactly}(1,\mathbf{Y_R}_j)$, i.e., exactly one variable within $\mathbf{Y_R}_j$ is \textit{True} ~\parencite{stevens2017teaching}. For example, in the case where there are only two independent Boolean variables ($\mathbf{Y_R}_j=(Y_1,Y_2)$) the constraint is equivalent to $Y_1 \veebar Y_2$. Note that, if the Boolean variables are transformed into binary variables, this is equivalent to a cardinality constraint $\sum_{i \in S_j} y_i = 1$ \parencite{ram1988algorithm}.
\end{itemize}

\subsubsection{External Variable Reformulation: Illustrative Example}

\begin{figure}
    \centering
    \includegraphics[width=0.70\textwidth]{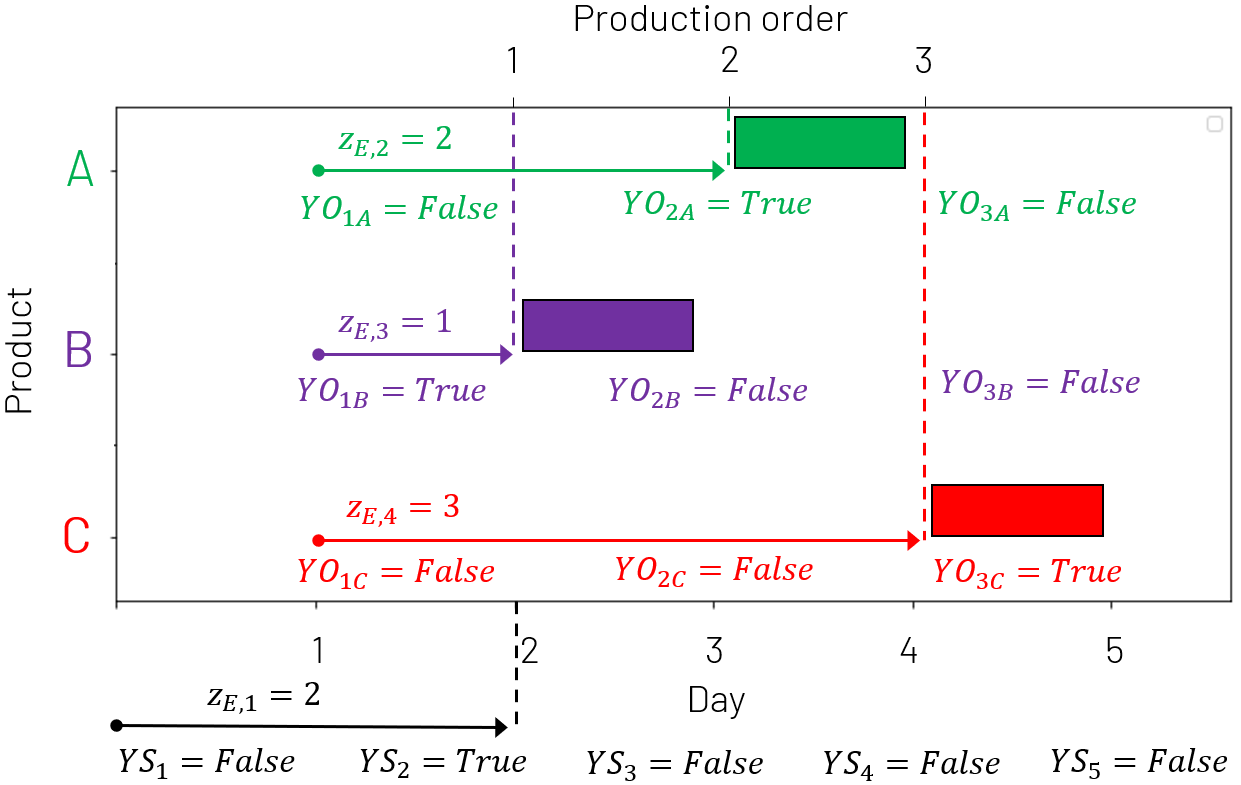}

    \caption{ Visualization of the external variable reformulation for an illustrative multi-product batch scheduling example. The figure explicitly displays the Boolean variables for starting time ($\mathbf{YS}$) and production order ($\mathbf{YO}_c$). Production begins on the second day, represented by $YS_2=True$ in black. The production order is B, A, and C, indicated by 
    $YO_{1B}=True$  (purple), $YO_{2A}=True$ (green), and $YO_{3C}=True$ (red). The maintenance variable $YM$ is not reformulated as it does not meet the necessary criteria for this transformation.}
    \label{fig:external_ref_example}
\end{figure}

To illustrate these requirements, consider the following example where a multi-product batch reactor that manufactures three products, $A$, $B$, and $C$. The goal is to determine an optimal starting date for each product within a five-day time horizon. Additionally, the production order for $A$, $B$, and $C$ must be established, subject to demand constraints. Another decision involves whether or not to perform routine maintenance before production begins.

To model this with Boolean variables, we define $YS_{t}, \forall \ t\in T=\{1,2,3,4,5\}$ to indicate the starting time, $YO_{pc}, \forall \ c\in C=\{A, B, C\}, \forall \ p\in P=\{1,2,3\}$ to represent the production order, and $YM$ to indicate the decision about performing maintenance.  The constraints of this problem dictate that there must be only one starting day
($\texttt{Exactly}(1,\mathbf{YS}=(YS_1, YS_2,..., YS_5))$) and that each product must be produced only once ($\texttt{Exactly}(1,\mathbf{YO}_c=(YO_{1c}, YO_{2c}, YO_{3c})) \ \forall \ c\in C$). 
These constraints imply that variables $\mathbf{YS},\mathbf{YO}_A,\mathbf{YO}_B$, and $\mathbf{YO}_C$ satisfy \textbf{Requirement 2}.
Moreover, since both $T$ and $P$ are ordered sets, we conclude the aforementioned group of variables also satisfies \textbf{Requirement 1}.
Hence, the vectors of independent Boolean variables can be grouped as $\mathbf{Y_R}=(\mathbf{YS},\mathbf{YO}_A,\mathbf{YO}_B,\mathbf{YO}_C)$, and reformulated with one external variable assigned to each vector in $\mathbf{Y_R}$. This means that ${z_{E,1}}$ is assigned to $\mathbf{YS}$, while ${z_{E,2}}, {z_{E,3}}$, and ${z_{E,4}}$ are assigned to $\mathbf{YO}_A,\mathbf{YO}_B$, and $\mathbf{YO}_C$, respectively.

The resulting reformulation is illustrated in Figure \ref{fig:external_ref_example} which shows the Boolean and external variable values of a possible solution. In this solution, production starts on Day 2, implying $YS_2 = True \Leftrightarrow  \mathbf{z}_{E,1}=2$, as indicated in black at the lower horizontal axis. 
Similarly, the upper horizontal axis indicates the production order where $B$ is produced first, followed by $A$, and then $C$.
This production order is represented by the reformulation $YO_{1B} = True \Leftrightarrow  \mathbf{z}_{E,2}=1, YO_{2A} = True \Leftrightarrow  \mathbf{z}_{E,3}=2$, and $YO_{3C} = True \Leftrightarrow  \mathbf{z}_{E,4}=3$. 
Finally, note that the Boolean variable $YM$ does not satisfy the requirements for reformulation, so it remains as is.

\subsubsection{External Variable Reformulation: Extension}

This work extends the external variable reformulation to general cases. To formally describe this approach, consider an optimization problem in the form of \eqref{eqn:gdp.general}.
If the GDP problem satisfies \textbf{Requirements 1} and \textbf{2} over the vectors of independent Boolean variables $\mathbf{Y_R}$, then one external variable can be assigned to each vector $\mathbf{Y_R}_j
, \forall \ j \in \{1,2,...,n_R\}
$.

\textbf{Requirement 1} indicates that each vector $\mathbf{Y_R}_j
, \forall \ j\in \{1,2,...,n_R\}$ 
 must be defined over a well-ordered set $S_j$. 
Not all Boolean variables defined over $S_j$ are required to belong to $\mathbf{Y_R}_j$. To account for this, we introduce the subset $S'_j \subseteq S_j$ for each $j \in {1,2,...,n_R}$, representing the ordered sets in which the Boolean variables $\mathbf{Y_R}_j$ are declared. 
The vector of independent variables is then represented as $\mathbf{Y_R}_j=(Y_{Rj,S'_j(1)},Y_{Rj,S'_j(2)},...,Y_{Rj,S'_j(\lvert S'_j \rvert)} )$ where each $Y_{Rj, S'_j(a)}$ is a Boolean variable defined at position $a$ of the well-ordered set $S'_j$ and included in vector $\mathbf{Y_R}_j$. 
It is important to note that  Boolean variables $Y_{Rj, S'_j(a)}
 , \forall \ j\in \{1,2,...,n_R\},
\ \forall \ a \in \{1,2,...,\lvert S'_j \rvert\}$ can be defined over other sets aside from $S'_j
, \forall \ j\in \{1,2,...,n_R\}
$. In other words, these variables may have multiple indices in the algebraic model formulation.

\textbf{Requirement 2} indicates that $\Omega(\mathbf{Y}) = True$ in \eqref{eqn:gdp.general} must contain partitioning constraints of the form $\texttt{Exactly}(1,\mathbf{Y_R}_j)
, \forall \ j\in \{1,2,...,n_R\}
$. 
Combining both requirements allows to define Boolean variables in $\mathbf{Y_R}$ as a function of $n_R$ external variables $z_{E,j}
 ,\forall \ j\in \{1,2,...,n_R\}
$ as,
\begin{equation}\label{eqn:reformulation}
    \begin{aligned}
        Y_{Rj,S'_j(a)} \iff z_{E,j}=a \quad \forall \ j \in \{1,2,...,n_R\}, \forall \ a \in \{1,2,...,\lvert S'_j \rvert\}
    \end{aligned}
\end{equation}

effectively expressing the external variables $z_{E,j}$ based on the values of Boolean variables $Y_{Rj, S'_j(a)}$.
From this reformulation, the upper and lower bounds of the external variables can be directly inferred from the sets of ordered positions   $\{1,2,...,\lvert S'_j \rvert\}
, \forall \ j \in \{1,2,...,n_R\}
$. These bounds are defined as:
\begin{equation}\label{eqn:ref_bounds}
    \begin{aligned}
 1 	\leq z_{E,j} \leq \lvert S'_j \rvert \quad \forall \ j \in \{1,2,...,n_R\} 
    \end{aligned}
\end{equation}





The general external variable reformulation is given by equations (\ref{eqn:reformulation}) and (\ref{eqn:ref_bounds}).
Next, we proceed to derive a simpler reformulation that follows from the special case when all the disjunctions are defined over well-ordered sets.
First, note that \textbf{Requirement 2} is naturally satisfied by the disjunctions in a standard \eqref{eqn:gdp.general} formulation. 
This arises from the fact that the exclusivity requirement in disjunctions enforces constraints of the form $\texttt{Exactly}(1,\mathbf{Y}_k), \forall \ k\in K$, where vector $\mathbf{Y}_k$ contains Boolean terms $Y_{ik}, \forall \ i\in D_k$. 
Consequently, if each disjunction $k$ represented an ordered decision over the well-ordered set $D_k$, then \textbf{Requirement 1} would be directly satisfied, allowing to reformulate a standard \eqref{eqn:gdp.general} problem following the guidelines in equations (\ref{eqn:reformulation}) and (\ref{eqn:ref_bounds}) to instead obtain:
\begin{equation}\label{eqn:reformulationGDP}
    \begin{aligned}
        Y_{D_k(a),k} \iff z_{E,k}=a \quad \forall \ k \in K, \forall \ a \in \{1,2,...,\lvert D_k \rvert\}
    \end{aligned}
\end{equation}
\begin{equation}\label{eqn:ref_boundsGDP}
    \begin{aligned}
 1 	\leq z_{E,k} \leq \lvert D_k \rvert \quad \forall \ k \in K 
    \end{aligned}
\end{equation}
respectively. 
In this case, there are as many external variables as disjunctions $k\in K$ in the formulation, making index $j$ interchangeable with disjunction index $k$.
Similarly, ordered subsets $S'_j$ correspond to disjunct sets $D_k$. 
For this reason, indices $i$ in $Y_{ik}$ are replaced by ordered index $D_k(a)$ in equation \eqref{eqn:reformulationGDP}. 
In practice, not every Boolean variable in the formulation fulfills the requirements to be reformulated with external variables as suggested by equation \eqref{eqn:reformulationGDP}.
In addition, ordered discrete structures may appear outside disjunctions, e.g., within  $\Omega(\mathbf{Y}) = True$.
Therefore, the reformulation in equations (\ref{eqn:reformulation}) and (\ref{eqn:ref_bounds}) is more general and practical than equations (\ref{eqn:reformulationGDP}) and (\ref{eqn:ref_boundsGDP}).    

Generalizing the reformulation established in previous research,
we adapt the proposed reformulation to handle variables $\mathbf{Y_R}_j$ defined over ordered, but unevenly spaced sets $\forall \ j\in \{1,2,...,n_R\}$.
Instead of defining external variables based on the \textit{elements} of these ordered sets $S_j$ as in earlier works \cite{linan2020optimalI}, we propose defining external variables based on the \textit{positions} in the ordered sets $S'_j$. For each $j$, the set of positions is denoted as ${1, 2, \dots, \lvert S'_j \rvert}$, where the distance between consecutive elements is equal to one.
This change avoids potential issues with solutions defined by isolated discrete elements, since the nearest-neighborhood search in the LD-SDA is now defined over positions instead of elements.

To illustrate the problem that may arise when considering the reformulation in terms of elements, consider an unevenly spaced set $S'_1=\{0,1,2,7,10\}$, and its corresponding Boolean variables $\mathbf{Y_R}_1=(Y_{R1,0}, Y_{R1,1}, Y_{R1,2}, Y_{R1,7}, Y_{R1,10})$, with the partitioning constraint $\texttt{Exactly}(1,\mathbf{Y_R}_1)=True$.
Suppose the incumbent point is $z_{E,1}=4$, which corresponds to $Y_{R1,7}=\textit{True}$. A neighborhood search around the incumbent would explore $z_{E,1}=3$ and $z_{E,1}=5$.
According to previous definitions \parencite{linan2020optimalI}, this search would attempt to set
$Y_{R1,3}=\textit{True}$ or $Y_{R1,5}=\textit{True}$, while assigning the remaining Boolean variables in $\mathbf{Y_R}_1$ to \textit{False}.
Since $3 \notin S'_1$ and $5 \notin S'_1$, both neighboring points would be treated as infeasible, causing the search to stop prematurely, identifying $z{E,1} = 4$ as a local optimum.
Our proposed definition resolves this issue by conducting the reformulation over set positions. In this case, the neighboring points $z_{E,1}=3$ and $z_{E,1}=5$ would correspond to positions over $\mathbf{Y_R}_1$, meaning that either $Y_{R1,2}$ or $Y_{R1,10}$ would be set to \textit{True}, instead of $Y_{R1,3}$ or $Y_{R1,5}$. This allows the discrete search to continue without prematurely declaring a local optimum.

\subsection{GDP Decomposition Using External Variables}
\label{sec:gdp_decomposition}
The reformulation presented in the previous section allows to express some of the Boolean variables in the problem in terms of the external variables as $\mathbf{Y_R}=\mathbf{Y_R}(\mathbf{z_E})$, where $\mathbf{z_E}$ is a vector of $n_R$ external variables. 
The core idea of the LD-SDA is to move these external variables to an upper-level problem \eqref{eqn:upper} and the rest of the variables to a subproblem \eqref{eqn:lower}.
This decomposition allows taking advantage of the special ordered structure of the \textit{external variables} by using a Discrete-Steepest Descent Algorithm in the upper-level problem to explore their domain as explained in \S \ref{sec:algorithm}.
Once an external variable configuration is determined by D-SDA, a subproblem is obtained by only considering the active disjuncts of that specific $\mathbf{z_E}$ configuration. 
The formal definition of both problems is given as:
\begin{equation}\label{eqn:upper}
\begin{aligned}
\min_{\mathbf{z_E}} &\ f_{sub}(\mathbf{z_E}) \\
\text{s.t.} &\ \mathbf{z_E} \in \{ \underline{\mathbf{z_E}},\dots,\overline{\mathbf{z_E}} \} \subseteq \mathbb{Z}^{n_{z_{E}}} &\ (\text{From Eq.} \ \eqref{eqn:ref_bounds}) \\
\end{aligned}\tag{Upper}
\end{equation}

\begin{equation}\label{eqn:lower}
\begin{aligned}
s(\mathbf{{z_{E}}})=
\begin{cases}
f_{sub}(\mathbf{z_E})= &\ \min_{\mathbf{x},\mathbf{Y_N},\mathbf{z}} f(\mathbf{x},\mathbf{z}) \\ 
&\ \text{s.t.} \ \mathbf{Y_R}=\mathbf{Y_R}(\mathbf{z_E}) \ \ \ \ \ (\text{Fixed as shown in Eq.} \ \eqref{eqn:reformulation}) \\
&\ \ \ \ \ \ \ \ \mathbf{g}(\mathbf{x}, \mathbf{z}) \leq \mathbf{0}\\
&\ \ \ \ \ \ \ \ \Omega(\mathbf{Y}) = True \\
&\ \ \ \ \ \ \ \ \bigveebar_{i\in D_k} \left[
    \begin{gathered}
    Y_{ik} \\
    \mathbf{h}_{ik}(\mathbf{x}, \mathbf{z})\leq \mathbf{0}\\
    \end{gathered}
\right] \quad k \in K\\
&\ \ \ \ \ \ \ \ \mathbf{x}  \in [\underline{\mathbf{x}},\overline{\mathbf{x}}] \subseteq \mathbb{R}^{n_x}; \mathbf{Y} \in \{False,True\}^{n_y}; \mathbf{z} \in \{ \underline{\mathbf{z}},\dots,\overline{\mathbf{z}} \} \subseteq \mathbb{Z}^{n_z}\\

\end{cases}
\end{aligned}\tag{Lower}
\end{equation}

In problem \eqref{eqn:upper}, a value for the objective function $f_{sub}(\mathbf{z_E})$ is obtained by the optimization of subproblem \eqref{eqn:lower}. 
Thus, $f_{sub}(\mathbf{z_E})$ is defined as an optimal objective function value found by optimizing the subproblem $s(\mathbf{z_E})$, obtained by fixing external variables fixed at $\mathbf{z_E}$.
If the subproblem \eqref{eqn:lower} is infeasible, $f_{sub}(\mathbf{z_E})$ is set to positive infinity by convention. Notably, the subproblems are reduced in the sense that they only consider the relevant constraints for the relevant external variable configuration $\mathbf{z_E}$.

A novel feature of the LD-SDA is its ability to handle various types of subproblems, extending previous versions \parencite{linan2020optimalI} that solely supported NLP subproblems.
In its most general form, the lower-level problem \eqref{eqn:lower} is a \eqref{eqn:gdp.general} with continuous ($\mathbf{x}$), discrete ($\mathbf{z}$) and non-reformulated Boolean ($\mathbf{Y_N}$) variables.
Consider the scenario where every Boolean variable can be reformulated (e.g., as shown in equations (\ref{eqn:reformulationGDP}) and (\ref{eqn:ref_boundsGDP})) or every non-reformulated variable $\mathbf{Y_N}$ is equivalently expressed in terms of $\mathbf{Y_R}$. Note that the latter situation may occur if all Boolean variables $\mathbf{Y_N}$ are determined within the subproblem upon fixing $\mathbf{Y_R}$, implying that logic constraints $\Omega(\mathbf{Y}) = True$ establish $\mathbf{Y_N}$ as functions of $\mathbf{Y_R}$.
In such a scenario, the resulting subproblem becomes an \eqref{eqn:minlp.general} with continuous ($\mathbf{x}$) and discrete variables ($\mathbf{z}$), or an NLP if there are no discrete variables ($\mathbf{z}$) in the formulation. 
In the following subsection, we introduce the LD-SDA as a decomposition algorithm that leverages the external variable reformulation and bi-level structure depicted so far.

\subsection{Logic-Based Discrete-Steepest Descent Algorithm}\label{sec:algorithm}

The Logic-Based Discrete-Steepest Descent Algorithm, as described in Algorithm \ref{alg:D-SDA}, solves a series of subproblems \eqref{eqn:lower} until a stopping criterion is satisfied. 
The LD-SDA can only start once the external variable reformulation of the problem has been performed. 
The external variables $\mathbf{{z_{E}}}$ are handled in an upper optimization level where the algorithm is performed.
To initialize, this method requires an initial fixed value of external variables $\mathbf{{z_{E,0}}}$, the value of the variables of its corresponding feasible solution ($ \mathbf{{x}_{0}}, \mathbf{{Y}_{0}},
\mathbf{{z}_{0}}$), and its respective objective function value $f_{sub}(\mathbf{{z_{E,0}}})$. 
Finding a starting feasible solution is beyond the scope of this work; however, it would be enough to have $\mathbf{{z_{E,0}}}$ and solve the subproblem $s(\mathbf{{z_{E,0}}})$ to find the rest of the required initial solution.
Note that problem-specific initialization strategies have been suggested in the literature, e.g., see \cite{linan2022optimal}.

The LD-SDA explores a neighborhood within the external variable domain; hence, the user must determine the type of neighborhood $k$ that will be studied.
As mentioned in \S\ref{sec:discrete_convex}, we only consider $k \in \{2, \infty\}$ (see Figure \ref{fig:neighborhoods}); nevertheless, other types of discrete neighborhoods can be considered \parencite{murota1998discrete}. Once $k$ has been selected, the neighborhood $N_k$ of the current point $\mathbf{{z_{E}}}$ is defined as $N_k(\mathbf{{z}_E}) = \{\mathbf{\alpha} \in \mathbb{Z}^{n_{z_{E}}}: \; 
\lVert \mathbf{\alpha} - \mathbf{{z}_E} \rVert_k \leq 1\}$ as given by Equation \ref{eq:neighborhood}. 
Similarly, the set of directions $\Delta_k$ from the point $\mathbf{{z_{E}}}$ to each neighbor $\mathbf{\alpha}$ is be calculated as $\Delta_k(\mathbf{{z}_E}) = \{\mathbf{d} : \; \mathbf{\alpha} - \mathbf{{z}_E} = \mathbf{d}, \forall \ \mathbf{\alpha} \in N_k(\mathbf{{z}_E})\}$ as dicatated by Equation \ref{eq:directions}.
 
The next step is to perform Neighborhood Search (see Algorithm \ref{alg:Neighborhood Search}), that consists of a local search within the defined neighborhood. 
Essentially, this algorithm solves $s(\mathbf{\alpha}) \; \forall \; \mathbf{\alpha} \in N_k(\mathbf{{z}_E})$ and compares the solutions found with the best incumbent solution $f_{sub}(\mathbf{z_E})$. 
If, in a minimization problem, $f_{sub}(\mathbf{z_E}) \leq f_{sub}(\mathbf{\alpha})  \; \forall \; \mathbf{\alpha} \in N_k(\mathbf{{z}_E})$ then, the current solution in $\mathbf{{z}_E}$ is a discrete local minimum (\emph{i-local} or \emph{s-local} depending on the value of \(k\)); otherwise, the steepest descent direction $\mathbf{d^*} = \mathbf{\alpha^*} - \mathbf{{z}_E}$ is computed, the algorithm moves to the best neighbor by letting $\mathbf{{z}_E} = \mathbf{\alpha^*}$ and performs a Line Search in direction $\mathbf{d^*}$. 
Note that for a neighbor to be considered the \emph{best neighbor} $\mathbf{\alpha^*}$, it must have a feasible subproblem, and a strictly better objective than both the incumbent solution and its corresponding neighborhood. 
 
The Line Search (see Algorithm \ref{alg:Line Search}) determines a point in the direction of steepest descent $\mathbf{\beta} = \mathbf{{z}_E} + \mathbf{d^*}$ and evaluates it. 
If the subproblem $s(\mathbf{\beta})$ is feasible and $f_{sub}(\mathbf{\beta}) < f_{sub}(\mathbf{z_E})$ then, let $\mathbf{{z}_E} = \mathbf{\beta}$ and perform the Line Search again until the search is unable to find a better feasible solution in direction $\mathbf{d^*}$. 
Once this occurs, the general algorithm should return to calculate $N_k(\mathbf{{z}_E})$ and $\Delta_k(\mathbf{{z}_E})$ to perform the Neighborhood Search again in a new iteration.
 
The LD-SDA will terminate once the Neighborhood Search is unable to find a neighbor $\mathbf{\alpha}$ with a feasible subproblem $s(\mathbf{\alpha})$ that strictly improves the incumbent solution as $f_{sub}(\mathbf{\alpha}) < f_{sub}(\mathbf{z_E})  \; \forall \; \mathbf{\alpha} \in N_k(\mathbf{{z}_E})$.
In that case, the point is considered a discrete \emph{i-local} or \emph{s-local} minimum, and the algorithm will return the values of both variables ($\mathbf{{x}},\mathbf{{Y}},\mathbf{{z}}, \mathbf{{z_E}}$) and the objective function $f_{sub}^*$ of the solution found.
The stopping criterion employed indicates that the current point has the best objective function amongst its immediate discrete neighborhood mapping  \parencite{murota1998discrete}.
In rigorous terms, the integrally local optimality condition can only be guaranteed after an $N_\infty$ exploration given that this is the neighborhood that considers the entire set of immediate neighbors (\emph{i-local} optimality).
Therefore, when using neighborhood $N_2$, it is up to the user to choose if the final solution $\mathbf{\hat{z_E}}$ is to be certified as integrally local by checking its $N_\infty(\mathbf{\hat{z_E}})$ neighborhood. 

\begin{figure}
    \centering
    \includegraphics[width=0.50\textwidth]{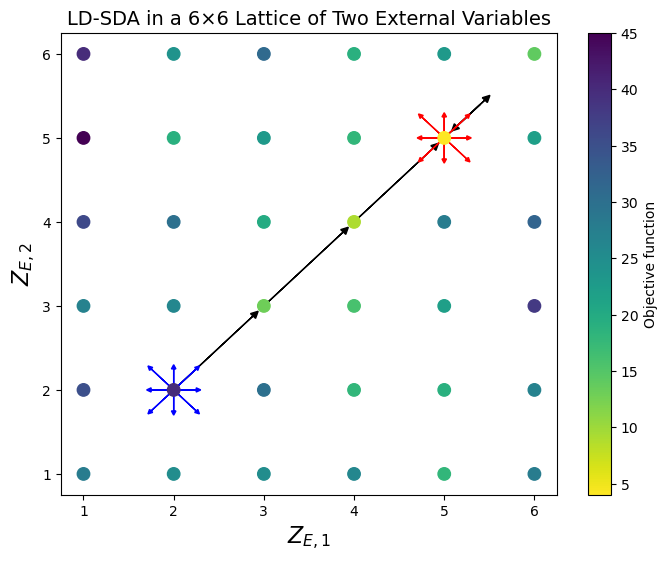}

    \caption{
    Visualization of the LD-SDA algorithm using $N_\infty$ in a two-variable discrete lattice example. In this example we initialize the algorithm begins at the initial point $(2,2)$. A Neighborhood Search within the neighborhood $N_\infty({(2,2)})$, represented by blue arrows, identifies the best neighbor as $(3,3)$, resulting in the steepest descent direction $\mathbf{d^*}={(1,1)}$. A Line Search, depicted with black arrows, follows this direction until reaching point $(5,5)$, where it stops given that $f_{sub}({(5,5)}) < f_{sub}({(5,5)}+\mathbf{d^*})$. A second Neighborhood Search in $N_\infty({(5,5)})$, shown with red arrows, determines $(5,5)$ is integrally local, terminating the algorithm.}
    \label{fig:6x6 feasible}
\end{figure}

Figure \ref{fig:6x6 feasible} illustrates and explains the LD-SDA executed in its entirety on a 6$\times$6 lattice of two external variables. For this example, the $\infty$-neighborhood is utilized and two different Neighborhood Searches are required.
Furthermore, a detailed pseudo-code for the LD-SDA is presented in the following section (see Algorithm \ref{alg:D-SDA}). 
Additional efficiency improvements and other implementation details are presented in \S \ref{sec:implementation details}.

\begin{algorithm}
\caption{Logic-Based Discrete-Steepest Descent Algorithm (LD-SDA)}
\label{alg:D-SDA}
\DontPrintSemicolon
  
  \KwInput{$k \in \{2, \infty\} $ ; 
   An external variable feasible solution
   $\mathbf{{z}_{E,0}}$}
  \KwData{Variable values associated with feasible solution $ \mathbf{{x}_{0}},
  \mathbf{{Y}_{0}},
   \mathbf{{z}_{0}}$}
   \tcc{Initialize}
  Set $\mathbf{{x}} \leftarrow \mathbf{{x}_{0}}; \;\;
		\mathbf{Y} \leftarrow \mathbf{{Y}_{0}}; \;\;
		\mathbf{z} \leftarrow \mathbf{{z}_{0}}; \;\;
		\mathbf{z_E} \leftarrow \mathbf{{z}_{E,0}}$ \\
  Solve subproblem: $f_{sub}^* \leftarrow f_{sub}(\mathbf{z_E})$ \\
  Set \texttt{neighborSearching} $\leftarrow True$ \\
  Generate initialization: $\gamma_{init} \leftarrow \{\mathbf{{x}},\mathbf{{Y}},\mathbf{{z}},\mathbf{{z_E}}\}$ \tcp*{Optional} 
  Initialize set of explored point in lattice as $G \leftarrow \{\mathbf{{z}_E} \}$\tcp*{Optional} \label{line:neighbor_search_check}
  \tcc{This cycle performs Neighborhood Search either when the algorithm starts (after initialization) or when Line Search does not improve the incumbent solution}
  \While {\normalfont{\texttt{neighborSearching} is} $True$}
  {     
  \tcc{Find the current neighborhood $N_k(\mathbf{{z}_E})$ and  directions $\Delta_k(\mathbf{{z}_E})$ to execute Neighborhood Search}
  Compute $N_k(\mathbf{{z}_E}) = \{\mathbf{\alpha} \in \mathbb{Z}^{n_{z_{E}}}: \; \lVert \mathbf{\alpha} - \mathbf{{z}_E}\rVert_k \leq 1\}$ \\
  Compute $\Delta_k(\mathbf{{z}_E}) = \{\mathbf{d} : \; \mathbf{\alpha} - \mathbf{{z}_E} = \mathbf{d}, \forall \ \mathbf{\alpha} \in N_k(\mathbf{{z}_E})\}$ \\
  
  \tcc{Perform the Neighborhood Search by evaluating an comparing every $f_{sub}(\mathbf{\alpha})$}
  
  $f_{sub}^* \;;\; \mathbf{{z_E}} \;;\; \mathbf{{d^*}} \;;\; \texttt{improvedDuringNS} \;;\; 
\gamma_{init} \leftarrow$ \textbf{Neighborhood Search}($f_{sub}^*, \mathbf{{z}_E}, N_k(\mathbf{{z}_E}), \Delta_k(\mathbf{{z}_E})$) 
  
  \tcc{Check for improvement during Neighborhood Search}
  
  \If{\normalfont{\texttt{improvedDuringNS} is} $True$}
    {
    \tcc{If so, perform Line Search in direction $\mathbf{d^*}$ until the incumbent does not improve}
         Set \texttt{lineSearching} $\leftarrow True$ \\
         \While {\normalfont{\texttt{lineSearching }is} $True$}
            {  
                $f_{sub}^* \;;\; \mathbf{{z_E}} \;;\; \texttt{improvedDuringLS} \;;\; 
\gamma_{init} \leftarrow$ \textbf{Line Search}($f_{sub}^*,\mathbf{{z}_E}, \mathbf{d^*}$)  \\
                \tcc{Check if the current solution was not improved during Line Search}
                \If{\normalfont{\texttt{improvedDuringLS} is} $False$}
                    {
                    \tcc{If so, stop Line Search}
                        Set \texttt{lineSearching} $\leftarrow False$ \\
                    }
             }
    }
    \Else
    {
        \tcc{If not, stop Neighborhood Search to terminate the algorithm and return the solution}
    	 Set \texttt{neighborSearching} $\leftarrow False$
    }
  }
   \KwOutput{$f_{sub}^* \;;\; \mathbf{{x}} \;;\; \mathbf{{Y}} \;;\; \mathbf{{z}} \;;\; \mathbf{{z_E}}$}

\end{algorithm}



\subsubsection{Neighborhood Search and Line Search}
\label{sec:neighbor_and_line_search}

In a general sense, the Neighborhood Search algorithm is a local search around the immediate neighborhood of discrete variables from a starting point $\mathbf{{z}_E}$. 
Therefore, the neighborhood $N_k(\mathbf{{z}_E})$ and the set of distances corresponding to each neighbor $\Delta_k(\mathbf{{z}_E})$ must be computed before starting the exploration. 
This algorithm solves the subproblems $s(\mathbf{\alpha}) \; \forall \ \mathbf{\alpha} \in N_k(\mathbf{{z}_E})$ and compares their objective function; if feasible, with the best incumbent solution found by Neighborhood Search $f_{sub}^{NS}$ so far.

The Neighborhood Search algorithm determines whether a new neighbor $\mathbf{\alpha}$ improves upon the current solution $\mathbf{z_E}$ based on two criteria, which can be evaluated in relative or absolute terms.
The first criterion employs a strict less than ($<$) comparison and is applied when no neighbor has yet outperformed the current solution. 
This ensures the algorithm avoids transitioning to a neighbor with an identical objective value, preventing cycling between points with the same objective function. 
Further details on the non-cycling properties of the LD-SDA are discussed in \S \ref{sec:properties dsda}.

Once a neighbor $\mathbf{\alpha}$ improves the incumbent solution, a second criterion (employing a less-than-or-equal-to ($\leq$) comparison) is utilized. This allows the algorithm to consider multiple neighbors $\mathbf{\alpha}$'s  with the same objective value. 
If more than one neighbor achieves the best solution, a tie-breaking strategy based on a maximum Euclidean distance lexicographic heuristic is used. The Euclidean distance is computed as $\texttt{dist} = \lVert \mathbf{\alpha} - \mathbf{{z}_E}\rVert_2\; \forall \ \mathbf{\alpha} \in N_k(\mathbf{{z}_E})$, favoring the first-found "most diagonal" path.
 These diagonal routes, which do not exist in $N_2$ neighborhoods, have proven effective in previous versions of the D-SDA \parencite{linan2020optimalI, linan2020optimalII, linan2021optimal}. 

\begin{algorithm}
\caption{Neighborhood Search}
\label{alg:Neighborhood Search}
\DontPrintSemicolon
  
  \KwInput{$f_{sub}^* \;;\; \mathbf{{z}_E} \;;\; N_k(\mathbf{{z}_E}) \;;\; \Delta_k(\mathbf{{z}_E})$}
  
  \tcc{Initialize}
  
  Relative tolerance $\epsilon 
  $ \\
    Set \texttt{improvedDuringNS} $\leftarrow False; \; \; 
    \mathbf{d^*} \leftarrow \mathbf{0}; \; \;
	f_{sub}^{NS} \leftarrow \infty$ \\
    Set $\texttt{dist}^* \leftarrow 0$ \tcp*{Optional}
    
\For{\normalfont{\textbf{every} $\mathbf{\alpha_i} \in N_k(\mathbf{{z}_E})$ }}
{
\tcc{Optional: Check if the neighbor was already evaluated in a previous iteration}
\If{$\mathbf{\alpha_i} \in G$ $\label{line: neighbor_search_check_starts}$}
{
Go to line 4 with $\mathbf{\alpha_{i+1}}$
}
\Else{Append $\mathbf{\alpha_i}$ to $G$ $\label{line: neighbor_search_check_ends}$}

\tcc{Optional: Check if the neighbor is within external variable domain}
\If{$\mathbf{\alpha_i} \notin \mathbf{Z_E} = \{ \mathbf{1},\dots,\overline{\mathbf{z_E}} \}$}
{
Go to line 4 with $\mathbf{\alpha_{i+1}}$
}

\tcc{Create fixed subproblem}

    Create subproblem and fix with external variables $s(\mathbf{\alpha_i})$  \\
	Initialize $s(\mathbf{\alpha_i})$ using $\gamma_{init}$ \tcp*{Optional} 
	\tcc{Optional: Check feasibility of fixed external variables in $s(\mathbf{\alpha_i})$ with FBBT}
	
	\If{\normalfont{FBBT of $s(\mathbf{\alpha_i})$} detects infeasibility}
{
Go to line 4 with $\mathbf{\alpha_{i+1}}$
}
\tcc{Solve subproblem}
	Solve $s(\mathbf{\alpha_i})$
  
  \If{\normalfont{$s(\mathbf{\alpha_i})$} is feasible}
    {
        Set $f_{sub}^{NS} \leftarrow f_{sub}(\mathbf{\alpha_i})$ \\
        
		Set $\texttt{dist}_i \leftarrow \lVert \mathbf{\alpha_i} - \mathbf{{z}_E}\rVert_2$ \tcp*{Optional}   
	
	\tcc{Check if the algorithm has already improved the starting solution to choose the corresponding minimum improvement criterion}
	
		\If{\normalfont{\texttt{improvedDuringNS }is} $False$}
		{
		
		\tcc{Check if minimum improvement criterion is satisfied}
		
		
		\If{$f_{sub}^{NS} < f_{sub}^* $ \normalfont{or} $(f_{sub}^* - f_{sub}^{NS})/(\lvert f_{sub}^* \rvert+10^{-10}) > \epsilon$}
		  {
		  \tcc{Update with new best solution}
		   Set $f_{sub}^* \leftarrow f_{sub}^{NS}$; $\mathbf{d^*} \leftarrow \Delta_k(\mathbf{{z}_E})_\mathbf{i}$;
		   $\mathbf{{z}_E}\leftarrow\mathbf{\alpha_i}$ \\
		  Set $\texttt{improvedDuringNS} \leftarrow True$ \\
		  Set $\texttt{dist}^* \leftarrow \texttt{dist}_i$ \tcp*{Optional}
		  
		  Generate initialization:  $\gamma_{init} \leftarrow \{\mathbf{{x}},\mathbf{{Y}},\mathbf{{z}},\mathbf{{z_E}}\}$ \tcp*{Optional} 
	   	  }
		}
    \Else
        {
    	\tcc{Check if minimum improvement criterion is satisfied. There is an additional condition that implements the maximum Euclidean distance heuristic}
    	
    	
    	\If{$(f_{sub}^{NS} \leq f_{sub}^* $ \normalfont{or} $(f_{sub}^* - f_{sub}^{NS})/( \lvert f_{sub}^* \rvert+10^{-10}) \geq \epsilon)$ \textcolor{black}{\normalfont{\textbf{and}}  $\texttt{dist}_i \geq \texttt{dist}^*$}}
    	{
    	\tcc{Update with new best solution}
		   Set $f_{sub}^* \leftarrow f_{sub}^{NS}$; $\mathbf{d^*} \leftarrow \Delta_k(\mathbf{{z}_E})_\mathbf{i}$;
		   $\mathbf{{z}_E}\leftarrow\mathbf{\alpha_i}$ \\
		  Set $\texttt{improvedDuringNS} \leftarrow True$ \\
		  Set $\texttt{dist}^* \leftarrow \texttt{dist}_i$ \tcp*{Optional}
		  
		  Generate initialization:  $\gamma_{init} \leftarrow \{\mathbf{{x}},\mathbf{{Y}},\mathbf{{z}},\mathbf{{z_E}}\}$ \tcp*{Optional} 
    	}
      }
    }
}
\KwOutput{$f_{sub}^* \;;\; \mathbf{{z_E}} \;;\; \mathbf{{d^*}} \;;\; \texttt{improvedDuringNS} \;;\; 
\gamma_{init}$}
\end{algorithm}

The Line Search algorithm is a search in the steepest descent direction, determined by the direction of the best neighbor $\mathbf{d^*} = \mathbf{\alpha^*} - \mathbf{z_E}$.
This approach generates a point in the steepest descent direction $\mathbf{\beta} = \mathbf{z_E} + \mathbf{d^*}$ and solves the optimization subproblem $s(\mathbf{\beta})$ to obtain $f_{sub}^{LS}$. 
The algorithm moves to the point $\mathbf{\beta}$ if and only if, $s(\mathbf{\beta})$ is feasible and $f_{sub}^{LS} < f_{sub}^{*}$, adhering to the strict less than ($<$) improvement criterion. 
This criterion prevents revisiting previous points, thereby accelerating and ensuring convergence.
Again, more insights into the convergence properties of the LD-SDA are stated in \S \ref{sec:properties dsda}.
The Line Search process continues until there is no feasible point in the direction $\mathbf{d^*}$ that improves upon the incumbent solution.

\begin{algorithm}
\caption{Line Search}
\label{alg:Line Search}
\DontPrintSemicolon
  
  \KwInput{$f_{sub}^* \;;\; \mathbf{{z}_E} \;;\; \mathbf{d^*}$}
  
  \tcc{Initialize}
  
  Relative tolerance $\epsilon 
  $ \\
    Set \texttt{improvedDuringLS} $\leftarrow False; \; \; 
    \mathbf{\beta} \leftarrow \mathbf{{z}_E} + \mathbf{d^*}; \; \;
	f_{sub}^{LS} \leftarrow \infty$

\tcc{Optional: Check if the moved point $\beta$ was already evaluated in a previous iteration} $\label{line: line_search_check_starts}$
\If{$\mathbf{\beta} \in G$} 
{
Terminate algorithm
}
\Else{Append $\mathbf{\beta}$ to $G$ $\label{line: line_search_check_ends}$ } 

\tcc{Optional: Check if the moved point is within the external variable domain}
\If{$\mathbf{\beta} \notin \mathbf{Z_E} = \{ \mathbf{1},\dots,\overline{\mathbf{z_E}} \}$}
{
Terminate algorithm
}

\tcc{Create fixed subproblem}

    Create subproblem and fix with external variables  $s(\mathbf{\beta})$  \\
	Initialize $s(\mathbf{\beta})$ using $\gamma_{init}$ \tcp*{Optional} 
	\tcc{Optional: Check feasibility of fixed external variables in $s(\mathbf{\beta})$ with FBBT}
	
	\If{\normalfont{FBBT of $s(\mathbf{\beta})$} detects infeasibility}
{
Terminate algorithm
}
\tcc{Solve subproblem}
	Solve $s(\mathbf{\beta})$
  
  \If{\normalfont{$s(\mathbf{\beta})$} is feasible}
    {
        Set $f_{sub}^{LS} \leftarrow f_{sub}(\mathbf{\beta})$ \\
	
	\tcc{Check if minimum improvement criterion is satisfied}
	
  \If{$f_{sub}^{LS} < f_{sub}^* $ \normalfont{or} $(f_{sub}^* - f_{sub}^{LS})/(\lvert f_{sub}^* \rvert+10^{-10}) > \epsilon$}
		  {
		  \tcc{Update with new best solution}
		   Set $f_{sub}^* \leftarrow f_{sub}^{LS}$; 
		   $\mathbf{{z}_E}\leftarrow \mathbf{\beta}$ \\
		  Set $\texttt{improvedDuringLS} \leftarrow True$ \\
		  
		  Generate initialization:  $\gamma_{init} \leftarrow \{\mathbf{{x}},\mathbf{{Y}},\mathbf{{z}},\mathbf{{z_E}}\}$ \tcp*{Optional} 
	   	  }
	
    }

\KwOutput{$f_{sub}^* \;;\; \mathbf{{z_E}} \;;\; \texttt{improvedDuringLS} \;;\; 
\gamma_{init}$}
\end{algorithm}

\subsection{Logic-Based Discrete-Steepest Descent Algorithm Properties}
\label{sec:properties dsda}


The LD-SDA algorithm is guaranteed not to cycle, meaning it avoids re-evaluating the same solution candidates while searching for an optimal solution.
This is achieved by avoiding revisitation of previously solved subproblems, and carefully deciding when to move to the next incumbent solution.
Both Neighborhood Search (Algorithm \ref{alg:Neighborhood Search}) and Line Search (Algorithm \ref{alg:Line Search}) adhere to a minimum improvement criterion to ensure progress toward a better solution. 
As discussed in \S \ref{sec:neighbor_and_line_search}, this criterion is satisfied if and only if a strict less than ($<$) improvement is obtained.
Once this is met, the algorithms update the incumbent with the newly found solution, ensuring that only strictly better solutions are accepted.
An important aspect of this approach is that it excludes points with identical objective values as the incumbent, effectively preventing cycling between points with the same objective. 
By avoiding the reevaluation of points in the lattice, both Neighbor and Line Search avoid retracing steps, which guarantees convergence to a discrete local minimum while also saving computational time by not re-evaluating the same point.

The primary advantage of the LD-SDA over previous iterations of the D-SDA in how it leverages the structure of ordered Boolean variables for external variable reformulation, rather than ordered binary variables.
In the LD-SDA, the solution to the upper-level problem is the same as the one used in the D-SDA, involving a series of Neighbor and Line searches over the external variable lattice. 
However, the key distinction is that each lattice point in the LD-SDA upper-level problem corresponds to a reduced space GDP or (MI)NLP, obtained by fixing Boolean variables, and thereby disjunctions (see \S \ref{sec:gdp_decomposition}).
 This leads to a reduced subproblem that only considers relevant constraints, avoiding zero-flow issues and improving both numerical stability and computational efficiency.
 In contrast, previous versions of the D-SDA fixed binary variables to obtain NLP subproblems, that could potentially contain irrelevant constraints with respect to the current configuration of the Boolean variables yielding and ill-posed problem.
 Furthermore, additional algorithmic improvements with respect to previous versions of the D-SDA were added to the LD-SDA as discussed in 
\S \ref{sec:implementation details}.


\subsection{Equivalence to Other Generalized Disjunctive Programming Algorithms}


While LD-SDA exhibits different features compared to other GDP algorithms, certain aspects of it remain equivalent to them. 
Notably, akin to other logic-based approaches, LD-SDA addresses (MI)NLP subproblems containing only the constraints of active disjunctions, thereby excluding irrelevant nonlinear constraints. 
Each method employs a mechanism for selecting the subsequent (MI)NLP subproblem, typically based on a search procedure.
In LOA, this mechanism involves solving a MILP problem subsequent to reformulating Problem \eqref{eqn:MainLOA}.
On the other hand, LBB determines a sequence of branched disjunctions for each level \(l\), \(KB_l\), based on a predetermined rule known as \emph{branching rule}.
In contrast, LD-SDA utilizes Neighbor and Line Search algorithms to make this decision, solving Problem \eqref{eqn:upper} locally.

The LD-SDA employs an external variable reformulation to map Boolean variables into a lower-dimensional representation of discrete variables. 
While LD-SDA solves the upper-level problem through steepest descent optimization, this problem essentially constitutes a discrete optimization problem without access to the functional form of the objective.
Hence, in principle, this problem could be addressed using black-box optimization methods.

Transitioning from one point to another in the discrete external variable lattice involves changing the configuration of Boolean variables in the original problem, often modifying multiple Boolean variables simultaneously. 
Consequently, the LD-SDA can be viewed as a variant of LBB, where Neighbor and Line Searches act as sophisticated branching rules for obtaining Boolean configurations to fix and evaluate. 
Furthermore, improvements to this problem could be achieved by leveraging information from the original GDP problem. 
For instance, linear approximations of the nonlinear constraints of the GDP could be provided, although this would necessitate employing an MILP solver. 
By constructing such linearizations around the solutions of the subproblem \eqref{eqn:lower}, one could recover Problem \eqref{eqn:MainLOA} from LOA.

\section{Implementation Details}
\label{sec:implementation details}


The LD-SDA, as a solution method for GDP, was implemented in Python using Pyomo~\parencite{bynum2021pyomo} as an open-source algebraic modeling language.
Pyomo.GDP~\parencite{chenpyomo} was used to implement the GDP models and use their data structures for the LD-SDA.
The code implementation allows the automatic reformulation of the Boolean variables in the GDP into external variables and provides an efficient implementation of the search algorithms over the lattice of external variables.

\subsection{Automatic Reformulation}

In contrast to previous works for MINLP models~\parencite{linan2020optimalI}, the reformulation in (\ref{eqn:reformulation}) and (\ref{eqn:ref_bounds}) provides a generalized framework that is automated in the Python implementation developed in this work.
 Minimal user input is required for the reformulation process, with only the Boolean variables in $\mathbf{Y}$ defined over ordered sets $S_j, \forall \ j \in \{1,2,...,n_R\}$ needing specification. 
This reformulation allows fixing Boolean variables based on the values of \emph{external variables}. 
Moreover, additional Boolean variables can be fixed based on the values of the external variables, as users can specify those Boolean variables in $\mathbf{Y_N}$ that are equivalent to expressions of the independent Boolean variables $\mathbf{Y_R}$ through logic constraints $\Omega(\mathbf{Y}) = True$.

\subsection{Algorithmic Efficiency Improvements}

This section presents the four major efficiency improvements that are included in the algorithm and are indicated throughout the pseudo-codes in \S \ref{sec:algorithm} as \texttt{\textcolor{blue}{Optional}}. 

\subsubsection{Globally Visited Set Verification}
Due to the alternating dynamic between Line Search and Neighborhood Search, the LD-SDA often queues discrete points that were previously visited and evaluated. 
An example of this issue can be observed in Figure \ref{fig:6x6 feasible} where the second Neighborhood Search in $N_\infty({(5,5)})$, depicted in red, visits points ${(4,4)}$ and ${(6,6)}$ that had already been evaluated during Line Search (shown in black). 
The number of re-evaluated points depends on how close to the Neighborhood Search the Line Search stops, increasing proportionally with the number of external variables. 

Although re-evaluating points does not affect the convergence of the algorithm as discussed in \S \ref{sec:properties dsda}, it results in unnecessary additional computation that can be avoided.
This redundant evaluation existed in the previous versions of the D-SDA \parencite{ linan2020optimalI, linan2020optimalII, linan2021optimal} and can be rectified by maintaining a globally visited set  $G$ (line \ref{line:neighbor_search_check} of Algorithm \ref{alg:D-SDA}). 
Furthermore, before solving the optimization model for a particular point $\mathbf{\alpha}$ 
(lines \ref{line: neighbor_search_check_starts} to \ref{line: neighbor_search_check_ends} of Algorithm \ref{alg:Neighborhood Search}) or $\mathbf{\beta}$ (lines  \ref{line: line_search_check_starts} to \ref{line: line_search_check_ends} of Algorithm \ref{alg:Line Search}), the algorithm verifies if the point has already been visited. 
If so, the algorithm disregards that point and either proceeds to the next $\mathbf{\alpha}$ in the Neighborhood Search or terminates the Line Search algorithm.

\subsubsection{External Variable Domain Verification}

All external variables must be defined over a constrained box $\mathbf{Z_E} = \{ \mathbf{1},\dots,\overline{\mathbf{z_E}} \}$ (as shown in Eq. \eqref{eqn:ref_bounds}) that depends on the problem. 
For superstructure problems, this domain is bounded by the size of the superstructure, such as the number of potential trays in a distillation column, or the maximum number of available parallel units in a process. 
Similarly, for scheduling problems, the external variable domain can be given by the scheduling horizon.

External variables with non-positive values or exceeding the potential size of the problem, resulting in a lack of physical sense, should not be considered in the explorations. 
To prevent unnecessary presolve computations, the algorithm verifies if the incumbent point ($\mathbf{\alpha}$ or $\mathbf{\beta}$) belongs to $\mathbf{Z_E}$  before solving the optimization model, effectively avoiding consideration of infeasible subproblems. 
If during Neighborhood Search  $\mathbf{\alpha} \notin \mathbf{Z_E}$, the neighbor $\mathbf{\alpha}$ can be ignored, and the algorithm proceeds to explore the next neighbor. 
Similarly, if $\mathbf{\beta} \notin \mathbf{Z_E}$ while performing the Line Search, the algorithm should return to $\mathbf{z_E} = \mathbf{\beta} - \mathbf{d^*}$ and terminate. 
Returning to the example shown in Figure \ref{fig:6x6 feasible}, note that, for instance, $N_\infty((1,6)) = \{(2,6),(2,5),(1,5)\}$ given that points $\{(1,7),(2,7),(0,5),(0,6),(0,7)\}$ can be automatically discarded and considered infeasible since $\mathbf{z_{E,1}}, \mathbf{z_{E,2}} \in \{1,2,...,6\}$. 

\subsubsection{Fixed External Variable Feasibility Verification via FBBT}
The existence of external variables $\mathbf{z_E}$ within their respective bounds does not ensure feasibility in the subproblem $s(\mathbf{z_E})$. 
While external variables can encode a physical interpretation of the problem by representing specific positions within a well-ordered set, constraints concerning the rest of the problem must align with spatial information to achieve a feasible subproblem.
For instance, consider the distillation column (discussed in \S \ref{sec:column}) that has two external variables: one determining the reflux position $z_{E, R}$ and another determining the boil-up position $z_{E, B}$. 
The problem has an implicit positional constraint $z_{E, B} < z_{E, R}$,  indicating that the boil-up stage must be above the reflux stage when counting trays from top to bottom.

Throughout the algorithm, this type of discrete positional constraint, which relates external variables, is frequently violated when a particular  $\mathbf{z_E}$ is fixed in a subproblem $s(\mathbf{z_E})$.
This violation arises because these constraints are specified in the original GDP model in terms of Boolean variables. 
Consequently, after the external variable reformulation, fixed points in the discrete lattice may overlook the original logical constraints.

In previous works \parencite{linan2020optimalI,linan2020optimalII,linan2021optimal}, users were tasked with manually re-specifying these constraints in the domain of external variables.  However, this work aims to automate this requirement. Instead of solving infeasible models that consume computation time and may generate errors that terminate the algorithm, we used the Feasibility-Based Bound Tightening routine available in Pyomo. 
FBBT rapidly verifies feasibility over the fixed Boolean constraints, enabling the algorithm to identify subproblem infeasibility without executing a more resource-intensive MINLP or GDP presolve algorithm.
Now, if FBBT determines that a subproblem $s(\mathbf{z_E})$ is infeasible, the point $\mathbf{z_E}$ can be instantly disregarded.

\subsubsection{Re-Initialization Scheme}
The LD-SDA method incorporates an efficiency improvement that involves reinitializing from the best solution $\gamma_{init}=\{\mathbf{{x}},\mathbf{{Y}},\mathbf{{z}},\mathbf{{z_E}}\}$. Effective model initialization is crucial for achieving faster convergence, particularly as problems increase in size and complexity.
Initiating a discrete point with the solution of a neighboring point is intuitively reasonable. 
Since points in the external variable lattice are derived from Boolean configurations following ordered sets,  adjacent points are expected to yield very similar subproblems (e.g., adding an extra tray in a distillation column or starting a process one time step later). Therefore, initializing from an adjacent neighbor can offer an advantage of discrete-steepest descent optimization over black-box methods that search the lattice.

During Neighborhood Search with  $\mathbf{\alpha} \in N_k(\mathbf{{z}_E})$, all subproblems are initialized using the solved variable values from the best incumbent solution $s(\mathbf{z_E})$. 
Similarly, in Line Search each subproblem from the moved point $s(\mathbf{\beta})$ is initialized with the variable values of the best incumbent solution $s(\mathbf{z_E}= \mathbf{\beta} - \mathbf{d^*})$. 
This re-initialization methodology proved very efficient when integrated into the MINLP D-SDA in the rigorous design of a catalytic distillation column using a rate-based model \parencite{linan2021optimal}.

\section{Results}
\label{sec:results}

The LD-SDA is implemented as an open-source code using Python. 
The case studies, such as reactors, chemical batch processing, and binary distillation column design, are modeled using Python 3.7.7 and Pyomo 5.7.3 \parencite{bynum2021pyomo}.
The catalytic distillation column design case study was modeled using GAMS 36.2.0.
All the solvers used for the subproblems are available in that version of GAMS and were solved using a Linux cluster with 48 AMD EPYC 7643 2.3GHz CPU processors base clock frequency and 1.0 TB RAM. Although the Neighborhood Search can be trivially parallelized, this study limited experiments to a single thread. 
All the codes are available at \url{https://github.com/SECQUOIA/dsda-gdp}. 
The solvers used for the MINLP optimization are BARON \parencite{sahinidis1996baron}, SCIP 
 \parencite{achterberg2009scip}, ANTIGONE \parencite{misener2014antigone}, DICOPT \parencite{bernal2020improving}, SBB \parencite{bussieck2001sbb}, and KNITRO \parencite{byrd2006k}. 
KNITRO, BARON, CONOPT~\parencite{drud1994conopt}, and IPOPT~\parencite{wachter2006implementation} are used to solve the NLP problems.
The GDP reformulations and algorithms are implemented in GDPOpt \parencite{chenpyomo}.

\subsection{Series of Continuously Stirred Tank Reactors (CSTRs)}

        

\begin{figure}
    \centering
    \includegraphics[width=0.75\textwidth]{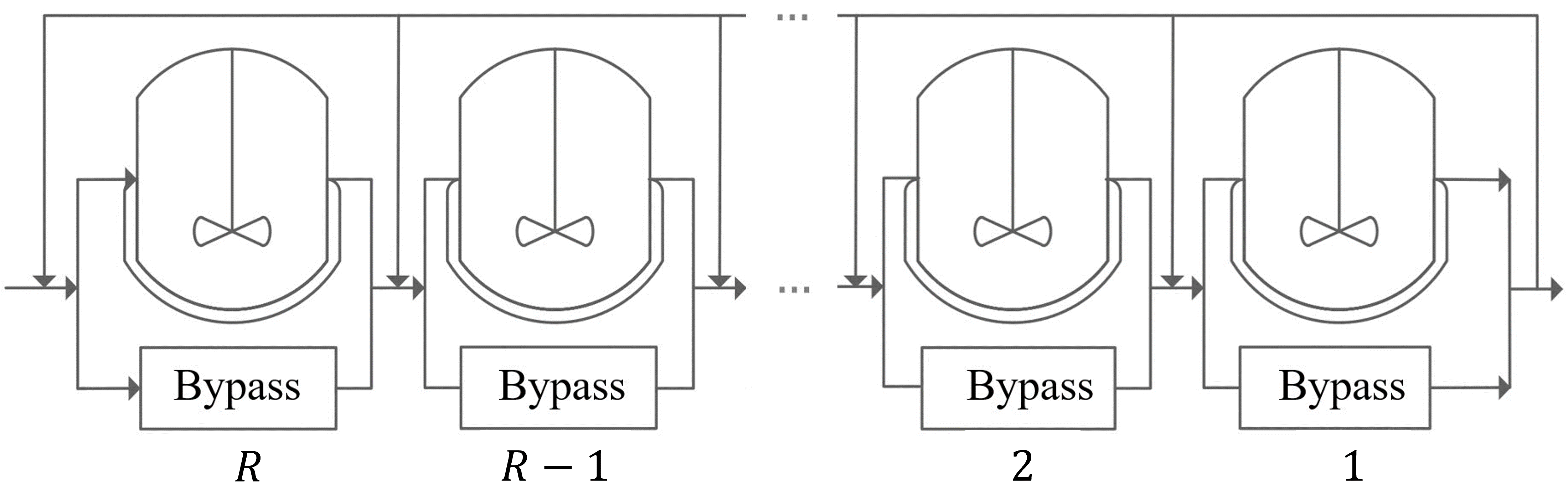}
    \caption{Visualization of a superstructure consisting of $R$ potential continuously stirred reactor tanks (CSTRs). The reactors are numbered starting from the product stream and counted in reverse. At each position, a reactor can either be present or replaced by a bypass. The configuration must be continuous, meaning no bypasses are allowed between two active reactors.}
    \label{fig:reactor_superstructure}
\end{figure}

Consider a reactor network adapted from  \parencite{linan2020optimalI}, consisting of a superstructure of \(R\) reactors in series (depicted in Figure  \ref{fig:reactor_superstructure}), where \(R\) represents the total number of potential reactors to install. 
The objective is to minimize the sum of reactor volumes.
The network involves an autocatalytic reaction \(A + B \rightarrow 2B\) with a first-order reaction rate, along with mass balances and reaction equations for each reactor. 
Logical constraints define the recycle flow location and the number of CSTRs in series to install. 
All installed reactors must have the same volume and a single recycle stream can feed any of them.
Interestingly, as the number of reactors increases and the recycle is placed in the first reactor, the system approximates to a plug-flow reactor, minimizing the total volume and providing an asymptotic analytical solution. 
We investigate this feature by varying the number of potential reactors \(R\).
For a detailed formulation of the reactor series superstructure, refer to Supplementary Material \ref{sec:cstr_formulation}.

\begin{figure}
    \centering
    \includegraphics[width=0.75\textwidth]{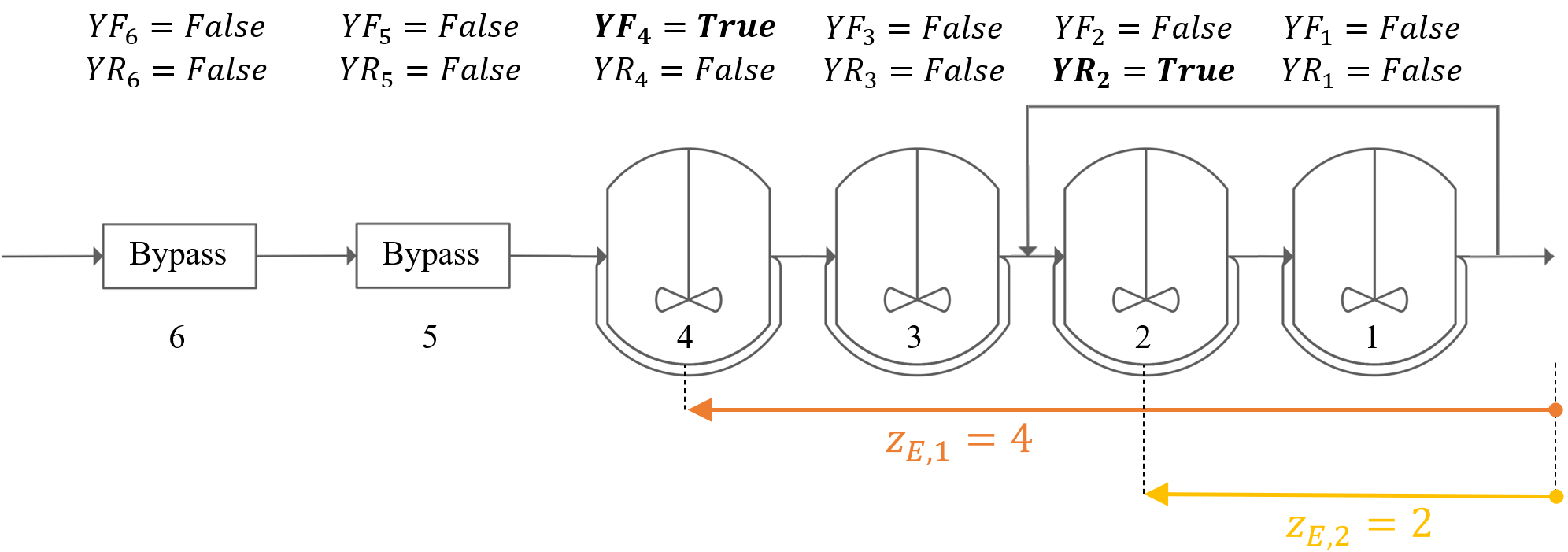}
    \caption{Visualization of a potential configuration for the CSTR superstructure with $R=6$. The figure explicitly displays the Boolean variables for the position of the feed ($\mathbf{YF}$) and recycle ($\mathbf{YR}$) streams. In this configuration, the feed enters at the fourth reactor, meaning $YF_4=True$, indicating the presence of four reactors in the superstructure. The recycle stream is positioned before the second reactor, meaning $YR_2=True$ .}
    \label{fig:reactor_reformulation}
\end{figure}

The external variables \(\mathbf{z_E} = (z_{E,1} : \text{No. of reactors (related with } \mathbf{YF}), z_{E,2} : \text{Recycle position (related with }  \mathbf{YR}))\) as shown in Figure \ref{fig:reactor_reformulation}  are the result of a complete reformulation of logic variables into external integer variables (detailed in Supplementary Material  \ref{app:cstr_reform}). 
This figure shows the binaries associated with the values of the ordered Boolean variables and their corresponding external variable mapping for an illustrative feasible solution, effectively indicating the reformulation $ YF_4 = True \Leftrightarrow  {z}_{E,1}=4$ and $YR_2 = True \Leftrightarrow  {z}_{E,2}=2$.

\begin{figure}
    \centering
    \includegraphics[width=0.70\textwidth]{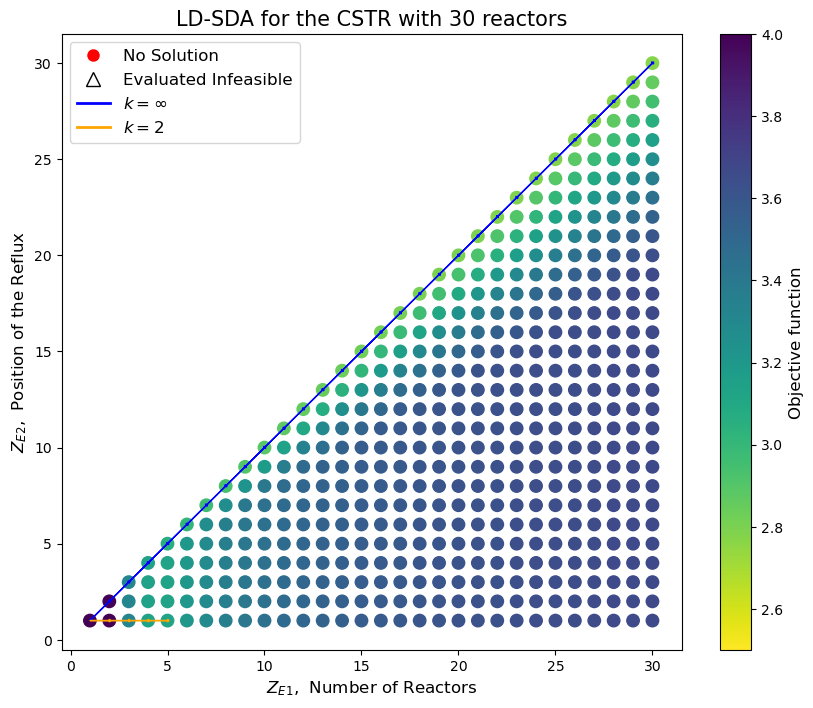}
    \caption{
    Visualization of the paths traversed by LD-SDA and the solutions found using both neighborhoods for a superstructure with $R=30$ CSTRs. LD-SDA with $N_2$ converged to $(5,1)$, a discrete $s$-local optimum in the lattice. In contrast, LD-SDA with $N_\infty$ continued to $(30,30)$, achieving the global optimal solution. This figure highlights the difference in convergence behavior between the two neighborhood strategies.}
    \label{fig:cstr_path}
\end{figure}

We analyze the paths and solutions generated by the LD-SDA under varying neighborhood selections, as depicted in Figure \ref{fig:cstr_path}.
Given \(R\) series of potential reactors, the problem is initialized with one reactor and its recycle flow. 
This initialization is represented in the integer variables lattice as a single reactor with a reflux position immediately behind it (\(\mathbf{z_E}=(1,1)\)).

For LD-SDA employing a \(k=2\) neighborhood search, the algorithm identifies \((2,1)\) is locally optimal and proceeds with the line search in the \(\mathbf{d^*}=(1,0)\) direction.
The algorithm continues the line search until \((5,1)\) as \((6,1)\) exhibits a worse objective. 
It searches among its neighbors, eventually moving and converging to a local optimal solution \((5,1)\).
In contrast, with LD-SDA utilizing a \(k=\infty\) neighborhood search, the algorithm finds that both \((2,1)\) and \((2,2)\)  yield the best solution within the first neighborhood explored. 
Employing the maximum Euclidean distance heuristic as a tie-break, this criterion selects \((2,2)\) as the new incumbent. 
Subsequently, a line search in the steepest direction \(\mathbf{d^*}=(1,1)\) proceeds until reaching \((R, R)\), representing the global optimal solution as it approximates the plug-flow reactor.

For a reactor series superstructure with \(R=30\), we performed external variable reformulation and fully enumerated the discrete points in a 30$\times$30 lattice.
Notably, only the lower-right triangle of the lattice is depicted, as points outside this region yield infeasible Boolean configurations. 
More specifically, these points indicate superstructures that have their recycle previous to an uninstalled reactor.
Local optimality was verified with both neighborhoods, revealing local minima for  \(k=2\) neighborhood at \((5,1)\), \((5,3)\) and \((r, r)\) $\forall \ $ \(r \in \{5, \dots, 29\}\), whereas the only locally optimal point for \(k=\infty\) was \((30,30)\).
Figure \ref{fig:cstr_path} illustrates the LD-SDA process for the 30 CSTR series, showing trajectories and local minimum points for all neighborhoods. 
The presence of multiple local optima with respect to both the $2$-neighborhood and the $\infty$-neighborhood suggests that this problem is neither separably convex, nor integrally convex.

\begin{figure}
    \centering
    \includegraphics[width=0.95\textwidth]{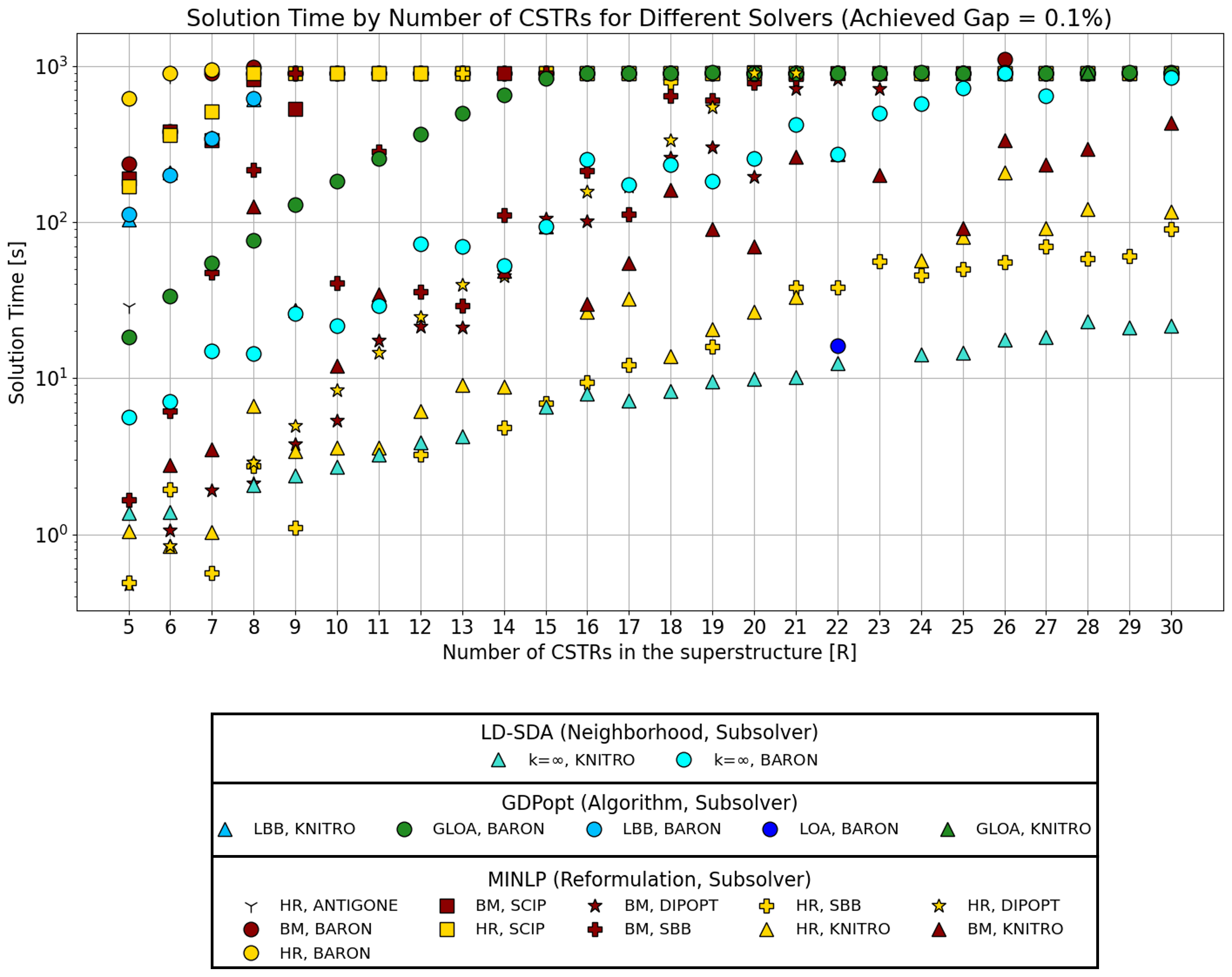}
    \caption{Computational solution times for different GDP solution strategies and solvers as the size of the superstructure $R$ increases. The figure includes only solutions that achieved a global optimum, which corresponds to the minimum total volume of the superstructure, analytically derived as the volume of a plug flow reactor as $R \rightarrow \infty$. It can be observed that LD-SDA with $N_{\infty}$ using KNITRO consistently reaches the global optimum the fastest for instances with more than 15 CSTRs.}
    \label{fig:reactor_0gap}
\end{figure}

For the CSTR series, various solver approaches are applied across different numbers of potential reactors (\(R\) ranging from 5 to 30). 
The solution approaches include MINLP reformulations, LBB, LOA, GLOA, and LD-SDA with two different neighborhoods.
Figure~\ref{fig:reactor_0gap} illustrates the comparison of solution times for each reactor superstructure size with different solvers. 

LD-SDA with \(k=2\) neighborhoods failed to achieve solutions within \(0.1 \% \) of the global optimum for any superstructure size, converging instead to the \((5,1)\) solution, as explained above. 
In contrast, methods that employed \(k= \infty\) attained the global minimum.
KNITRO was computationally more efficient than BARON when using \(k=\infty\), as BARON, being a global solver, incurred greater computational costs certifying global optimality for each NLP subproblem. 
Remarkably, even when using local solvers like KNITRO for the subproblem, \(k=\infty\) allowed LD-SDA to converge to global optimal solutions.

Among the logic-based methods in GDPopt, LBB for \(R=5\) achieved the global optimum. 
GLOA reached global optimal solutions up to \(R = 14\), but only when paired with the global NLP solver BARON.
Comparing MINLP reformulations, HR outperformed BM, with KNITRO being the most efficient subsolver.
While some MINLP reformulations exhibited faster solution times for smaller superstructures (up to nine reactors), LD-SDA, particularly with KNITRO, surpassed them for larger networks (from 15 reactors onwards). 
This trend suggests that LD-SDA methods are particularly well-suited for solving larger optimization problems, where solving reduced subproblems offers a significant advantage over monolithic GDP-MINLP approaches.


\begin{figure}
    \centering
    \includegraphics[width=1.00\textwidth]{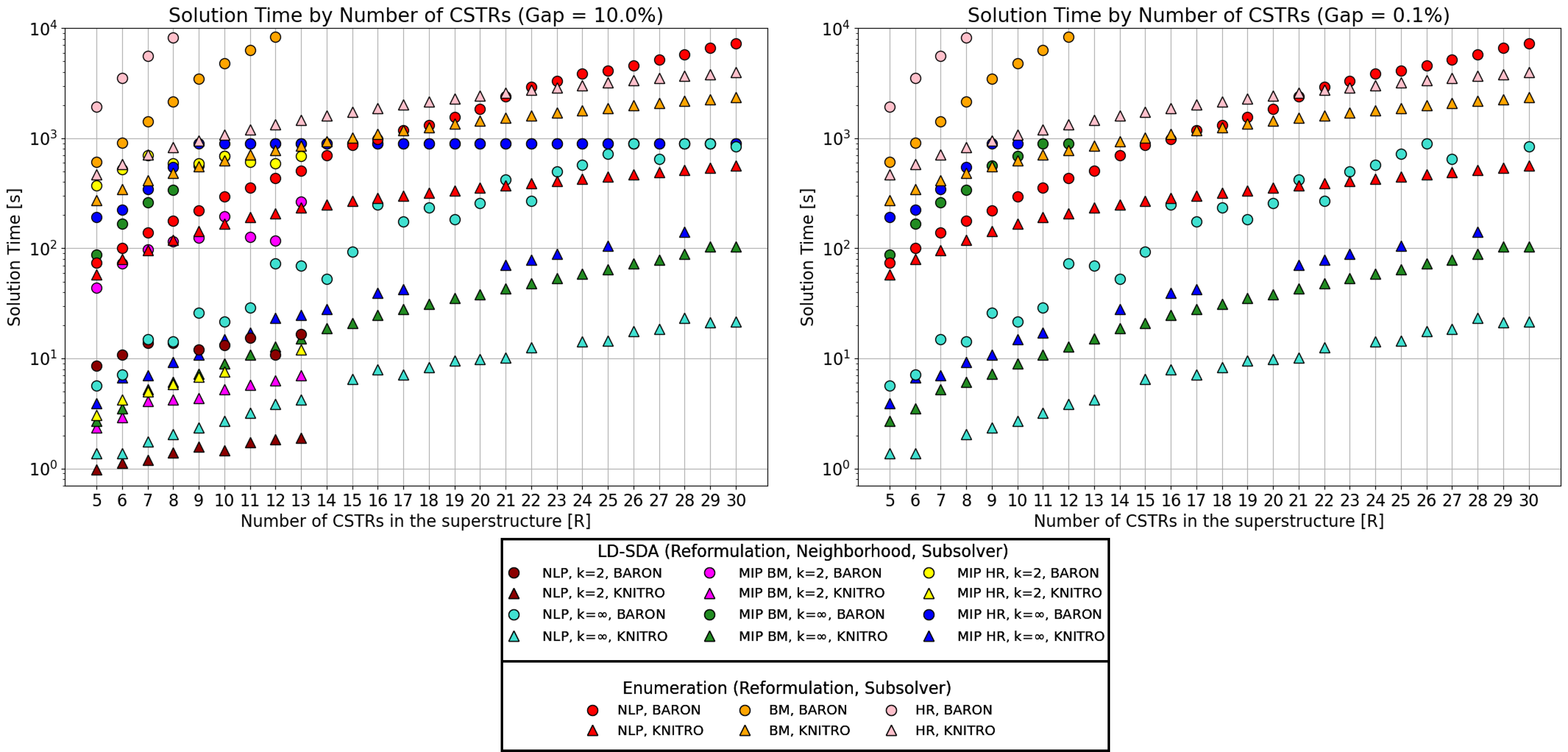}
    \caption{Computational solution times for different solution approaches derived from LD-SDA (including external variable enumeration) using various solvers as the size of the superstructure $R$ increases. The left subfigure shows methods that achieved a discrete local minimum with a 10\% optimality gap, while the right subfigure shows methods that reached an analytically proven global optimum. LD-SDA with $N_2$ and the KNITRO solver is the fastest method for smaller superstructures but fails to meet the 10\% optimality gap for structures larger than 13 CSTRs. For larger superstructures, LD-SDA with $N_{\infty}$ and KNITRO consistently achieves a global optimum and becomes the fastest approach among LD-SDA methods.}
    \label{fig:reactor_10gap}
\end{figure}

Figure \ref{fig:reactor_10gap} compares different algorithmic alternatives derived from LD-SDA. 
These include the algorithm discussed so far (referred in this example as NLP LD-SDA) where Boolean variables are fixed from external variables, leading to NLP subproblems considering only relevant constraints. 
 Another approach, which we refer to as MIP LD-SDA, is where inactive disjunctions are retained in subproblems, and mixed-binary reformulations (e.g., HR or BM) are applied to unresolved disjunctions, resulting in MINLP subproblems.
The third alternative is Enumeration, which involves reformulating external variables, fixing (or not) Boolean variables, and enumerating all lattice points instead of traversing them via steepest descent optimization.


LD-SDA and Enumeration methods exhibited faster performance when the mixed-binary reformulation was omitted. 
The inclusion of MIP transformations led to additional solution time, emphasizing the efficiency of solving GDP problems directly where reduced subproblems with solely relevant constraints are considered.
As anticipated, the Enumeration of external variables coupled with an efficent local solver like KNITRO achieved the global optimum. 
However, employing LD-SDA yielded the same result in significantly less time, showcasing the importance of navigating the lattice intelligently, like via discrete-steepest descent.

Among the LD-SDA approaches, \(2\)-neighborhood search methods were only effective when the tolerance gap between solutions was 
 \(10 \%\), while \(\infty\)-neighborhood search methods performed consistently across both gap thresholds.
The LD-SDA using \(k=2\) neighborhood converges to a local minimum that is more than \(0.1 \% \) away from the global optimal solution, and for larger instances, is beyond the \(10 \%\) optimality gap. 
Although LD-SDA with the \(k=\infty\) neighborhood search required more time compared to \(k=2\), it consistently converged to the global optimal point regardless of the superstructure size.


\subsection{Distillation Column Design for a Binary Mixture}
\label{sec:column}



We consider the single-unit operation design of an example distillation column in~\textcite{ghouse2018comparative}, that implements the simplified model provided by~\textcite{jackson2001disjunctive}. 
The objective is to design a distillation column to separate Toluene and Benzene while minimizing cost, which include both a fixed cost for tray installation and operational costs for the condenser and reboiler.
 
The column processes 100 gmol/s of an equimolar benzene-toluene mixture, aiming to achieve a minimum mole fraction of 0.95 for benzene in the distillate and 0.95 for toluene in the bottom product. 
To meet these requirements, the design must satisfy mass, equilibrium, summation, and heat (MESH) equations for each tray.
Each stage of the column is modeled using thermodynamic principles and vapor-liquid equilibrium, applying Raoult’s law and Antoine equation.

The continuous variables in the model include the flow rates of each component in both the liquid and vapor phases, the temperatures of each tray, the reflux and boil-up ratios, and the heat duties of the condenser and reboiler.
The logical variables account for the existence of trays and the positions of reflux and boil-up flows, with tray existence modeled using logical constraints related to these flow positions.
Previous studies from the literature  \parencite{ghouse2018comparative} set a maximum number of 17 potential trays and provide the initial position of the feed tray in the ninth stage (tray nine).

\begin{figure}
     \centering
     \begin{subfigure}[b]{0.3\textwidth}
         \centering
         \includegraphics[width=\textwidth]{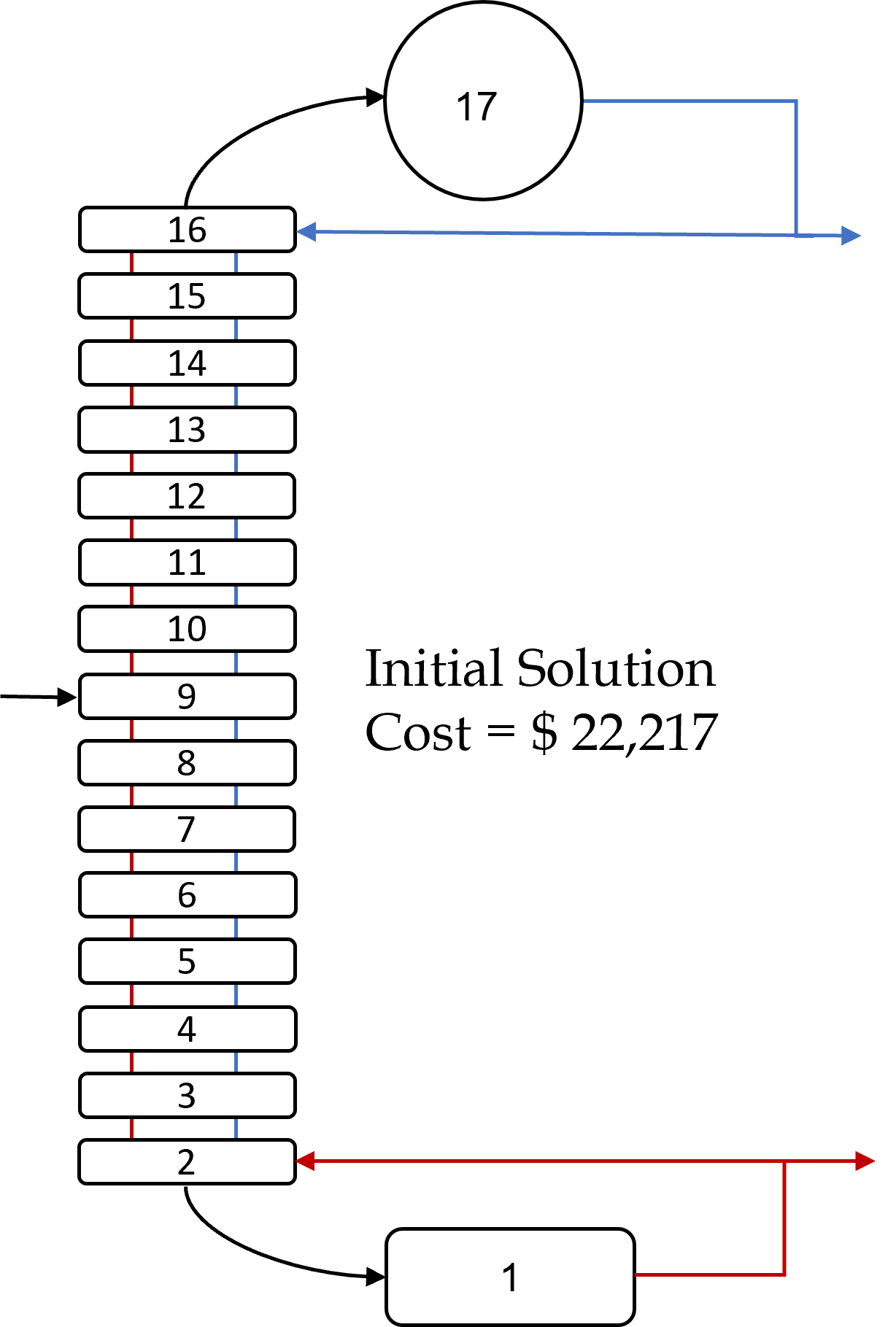}
         \caption{Distillation column configuration of the initial solution.}
         \label{fig: Initial Setting}
     \end{subfigure}
     \hfill
     \begin{subfigure}[b]{0.3\textwidth}
         \centering
         \includegraphics[width=\textwidth]{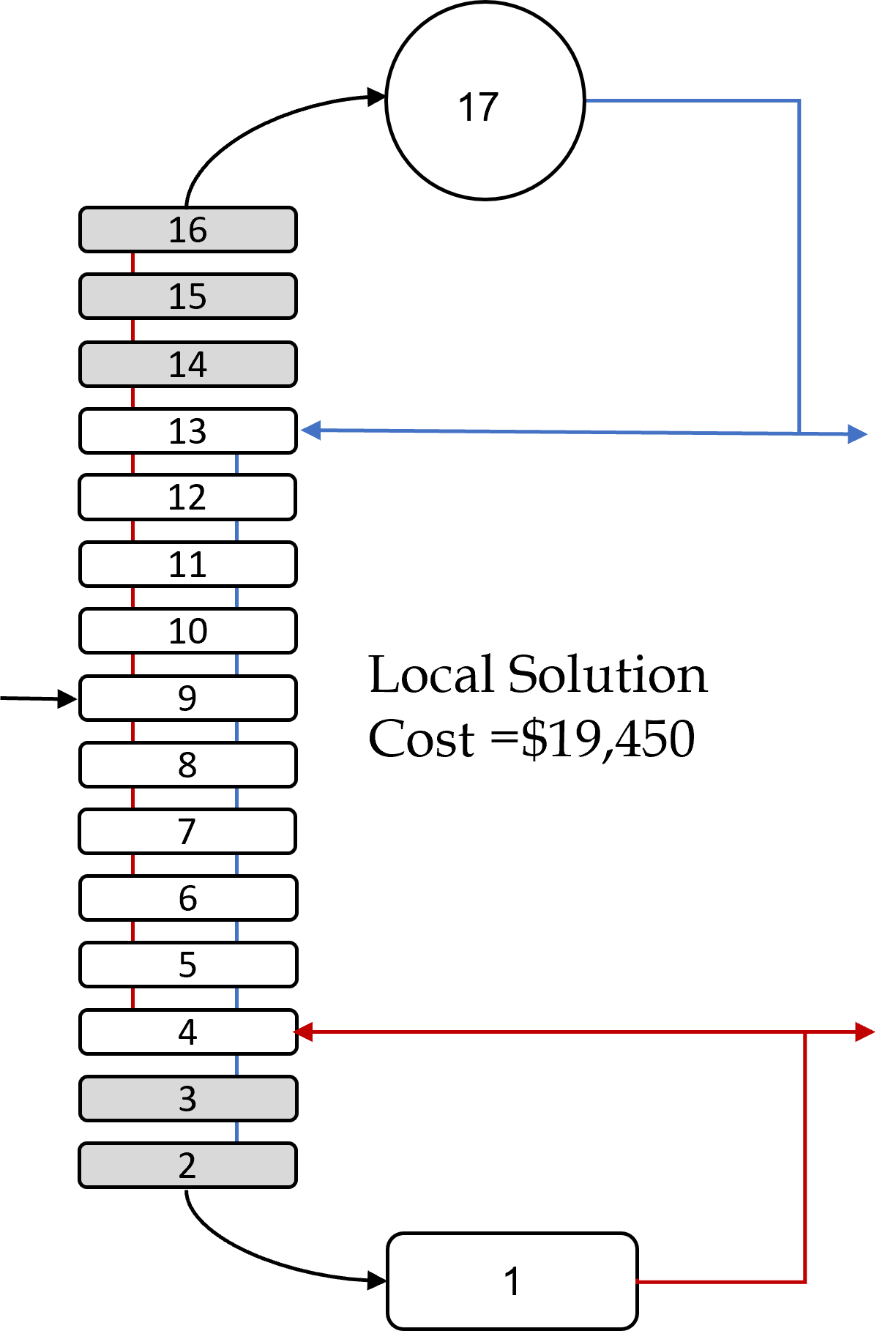}
         \caption{Distillation column configuration obtained with LD-SDA and $N_2$, which matches the solution reported in \parencite{ghouse2018comparative}.}
         \label{fig: Ghouse Solution}
     \end{subfigure}
     \hfill
     \begin{subfigure}[b]{0.3\textwidth}
         \centering
         \includegraphics[width=\textwidth]{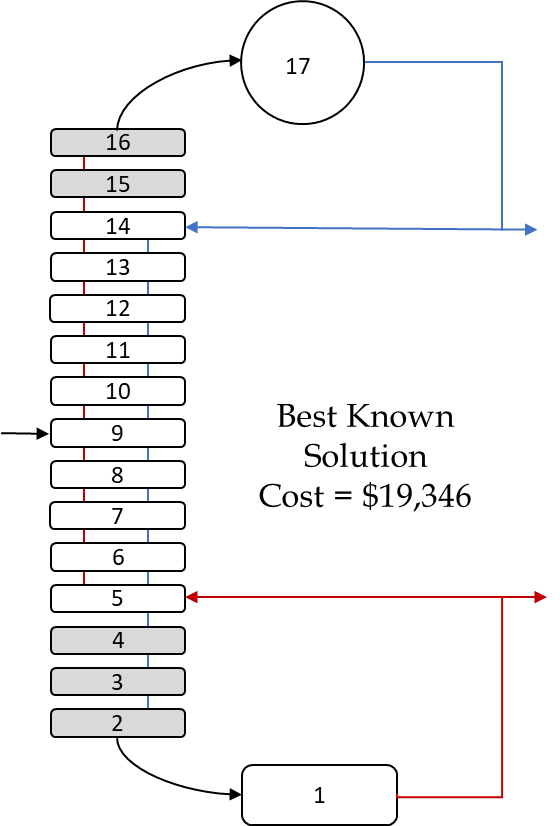}
         \caption{Distillation column configuration obtained with LD-SDA and $N_\infty$, which yields a new best solution for the problem.}
         \label{fig:LD-SDA Method}
     \end{subfigure}
        \caption{Visualization of the distillation columns obtained from the initial solution, the solution found using LD-SDA with $N_2$, and the solution using LD-SDA with $N_{\infty}$. The feed tray is fixed at stage nine, with existing trays displayed in white and bypassed trays shaded in gray. The boil-up position is highlighted in red, while the reflux position is marked in blue. Each subfigure shows its corresponding objective function.}
        \label{fig:Distillation Column }
\end{figure}

The distillation column optimization uses the LD-SDA method with different neighborhoods for the search. 
The problem is initialized with a column that has the reflux at tray \(16\) (top to bottom numbering, condenser being tray \(17\) and the reboiler being tray \(1\)) and the boil-up at the second tray, which we represent as \((15, 1)\).
The configuration for initialization is shown in Figure~\ref{fig: Initial Setting}, which corresponds to all possible trays being installed.

For the \(2\)-neighborhood, the LD-SDA converges to external variable configuration \((12, 3)\), with an objective value of \(\$ 19,450\) in only 6.3 seconds using KNITRO as subsolver, resulting in the design shown in Figure~\ref{fig: Ghouse Solution}.
Note this is the exact same solution reported by the GDP model from the literature \parencite{ghouse2018comparative}, that was solved using the LOA method. 

Regarding the \(\infty\)-neighborhood, the algorithm terminates at \((13, 4)\) with an objective of \(\$ 19,346\) after 8.6 seconds using KNITRO as subsolver, yielding the column design shown in Figure~\ref{fig:LD-SDA Method}.
In this case, the LD-SDA found the best-known solution to this problem, also found through a complete enumeration over the external variables, that took 42.7 seconds using KNITRO as a subsolver.
The same best-known solution could be found using GLOA with KNITRO as the NLP subsolver, but after 161.6 seconds.
In our results of the binary mixture distillation column design, we successfully identified a better solution by applying the LD-SDA to the GDP model, surpassing the optimal values previously reported in the literature.
This achievement highlights the efficacy of our approach, especially considering the limitations of the NLP formulation in guaranteeing global optimality, which we effectively navigated by employing the GDP framework.

\begin{figure}
    \centering
    \includegraphics[width=0.6\textwidth]{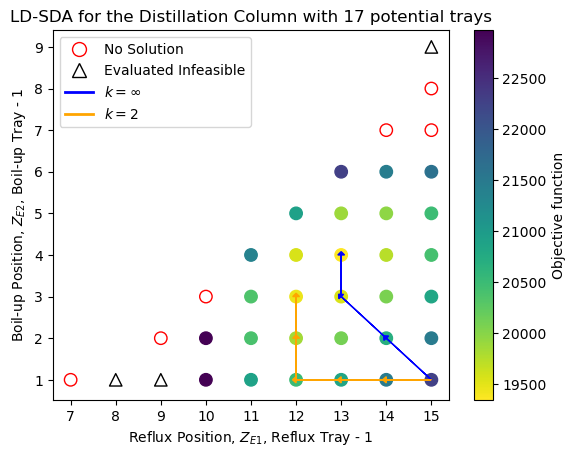}
    \caption{Visualization of the paths traversed by LD-SDA and the solutions found using both neighborhoods for the distillation column superstructure. 
    The lattice shows points proven to be infeasible with a global solver (black triangle) and points where no feasible solution was found before timing out (red circle). Both solutions were initialized with the largest column configuration. LD-SDA with $k=2$ converged to $(13,3)$, a discrete $s$-local optimal point in the lattice. Meanwhile, LD-SDA with $k=\infty$ advanced to $(13,4)$, yielding the best-known solution for the problem.}
    \label{fig:column_dsda}
\end{figure}

The trajectories traversed by the LD-SDA with both neighborhoods mentioned are depicted in Figure~\ref{fig:column_dsda}. 
Similarly, Table~\ref{tab:LDSDA_column} summarizes the previous design from the literature as well as the different columns obtained by the LD-SDA.

\begin{table}
    \centering
    \caption{Comparison of the optimal solutions found in the literature, using LD-SDA with $N_2$ and LDSDA with $N_{\infty}$. The LD-SDA with $N_2$ identified the same configuration as previously reported in the literature (with a minor numerical difference in the objective). In contrast, the LD-SDA with $N_{\infty}$ discovered a new best solution, effectively improving upon the existing results in the literature. This table illustrates the improvements gained by using more expansive neighborhoods in LD-SDA.}
\label{tab:LDSDA_column}
\begin{tabular}{cccc}
\toprule
\textbf{Solution Method}            & \textbf{LOA Ghouse et al. \parencite{ghouse2018comparative}} & \textbf{LD-SDA} $\mathbf{k=2}$ & \textbf{LD-SDA} $\mathbf{k= \infty}$ \\ \midrule
Objective [\$]              & 19,450                           & 19,449      & 19,346              \\ \midrule
Number of Trays             & 10                              & 10         & 10                 \\ \midrule
Feed Tray                   & 6                               & 6          & 5                  \\ \midrule
Reflux ratio                & 2.45                            & 2.45       & 2.01               \\ \midrule
Reboil ratio                & 2.39                            & 2.39       & 2.00               \\ \bottomrule      
\end{tabular}
\end{table}

\subsection{Catalytic Distillation Column Design}
Consider a catalytic distillation column design for the production of Ethyl tert-butyl ether (ETBE) from isobutene and ethanol.
In this work, two models are considered: one that uses equilibrium-based modeling in each of the separation and reactive stages, and another one that includes a rate-based description of the mass and energy transfer in all the stages~\parencite{linan2020optimalII}.
These models maximize an economic objective by determining the position of separation and catalytic stages along the column, together with a Langmuir-Hinshelwood-Hougen-Watson kinetic model for the chemical reaction, MESH equations for each one of the stages, and hydraulic constraints for the column operation.
The goal is to determine optimal operational variables such as reboiler and condenser heat duties and reflux ratio. 
Similarly, design variables such as column diameter, tray height, and downcomer specifications need to be defined. 
Finally, discrete design choices, meaning feed locations and positions of catalytic stages, must be selected.
A detailed description of the models is given in \parencite{linan2021optimal, bernal2018simultaneous}.

Previously in the literature, the economic annualized profit objective maximization of a catalytic distillation column to produce ETBE from butenes and ethanol was solved using a D-SDA \parencite{linan2021optimal}. 
Here, the authors demonstrated the difficulty of this design problem 
as several traditional optimization methods fail to obtain even a feasible solution \parencite{linan2020optimalI, linan2020optimalII}. 
In these papers, the D-SDA was used to solve the problem as an MINLP by fixing binary variables and including constraints of the form \(y_{ik}\mathbf{h}_{ik}(x) \leq 0\) to enforce the logic constraints. 
In this work, we demonstrate that approaching the problem disjunctively and employing LD-SDA leads to a faster solution of subproblems (as in Eq. \eqref{eqn:Sub}) as our method neglects the irrelevant and numerically challenging nonlinear constraints.

\begin{table}
\centering
\caption{Comparison of the optimal solutions and respective computational times using KNITRO for the catalytic distillation column and rate-based catalytic distillation column case studies, using D-SDA (from the literature) and LD-SDA with both $N_2$ and $N_{\infty}$. For the catalytic distillation column case, both neighborhoods in LD-SDA successfully found the same optimal solution as D-SDA but in one-third of the computational time. In the rate-based catalytic distillation column case, D-SDA was unable to find a solution, while LD-SDA, using both neighborhoods, found distinct optimal solutions.}
\label{tab:knitro_dantzig}
\begin{tabular}{ccccccccc} \toprule
                   & \multicolumn{4}{c}{ \textbf{Catalytic Distillation Column}}                & \multicolumn{4}{c}{\textbf{Rate-Based Catalytic Distillation Column}}     \\ \midrule
\textbf{Solution Method}    & \multicolumn{2}{c}{D-SDA: \parencite{linan2021optimal}} & \multicolumn{2}{c}{LD-SDA: This work} & \multicolumn{2}{c}{D-SDA: \parencite{linan2021optimal}} & \multicolumn{2}{c}{LD-SDA: This work} \\ \midrule
Neighborhood       & $k=2$   & $k=\infty$  &$k=2$   & $k=\infty$   & $k=2$  & $k=\infty$  & $k=2$    & $k=\infty$  \\ \midrule
Objective [\$/year] & 22,410   & 22,410       & 22,410   & 22,410        &   --   &    --    & 23,443.2      & 23,443.2       \\ \midrule
Time [s]           & 12.49    & 12.52       & 4.29      & 4.25        &  --  &   --     & 1089.31    & 1061.18       \\ \bottomrule
\end{tabular}
\end{table}

The models were implemented in GAMS.
Hence, for this problem, the reformulation and implementations of the algorithms were custom-made, as they did not rely on our implementation of LD-SDA in Python.
Given that only the relevant constraints were included for each problem, we could more efficiently obtain the same solution to each subproblem.
More specifically, as shown in Table~\ref{tab:knitro_dantzig}, the proposed LD-SDA method leads to speedups of up to three-fold in this problem when using KNITRO as a subsolver.
D-SDA was unable to even initialize the rate-based catalytic distillation column with KNITRO, while LD-SDA could find an optimal solution.
Moreover, note that the previous results using the D-SDA were already beating state-of-the-art MINLP solution methods, further demonstrating the advantages of the LD-SDA.

An important distinction for the LD-SDA is that it does not include all the constraints in each iteration, given that subproblems are reduced after the disjunctions are fixed.
This implies that not all variables are present in all iterations, preventing a complete variable initialization as the algorithm progresses.
These missing values for the variables might make converging these complex NLP problems challenging, which explains why the D-SDA and the LD-SDA sometimes yield different solutions.
Moreover, the solver KNITRO reported that the initial point was infeasible for the more complex NLP problem involving rate-based transfer equations.
Using that same initialization, yet using the Logic-Based D-SDA, the model could not only be started, but it converged to the same optimal solution reported in~\parencite{linan2020optimalII}.


\begin{figure}
     \centering
     \begin{subfigure}[b]{0.33\textwidth}
         \centering
         \includegraphics[width=\textwidth]{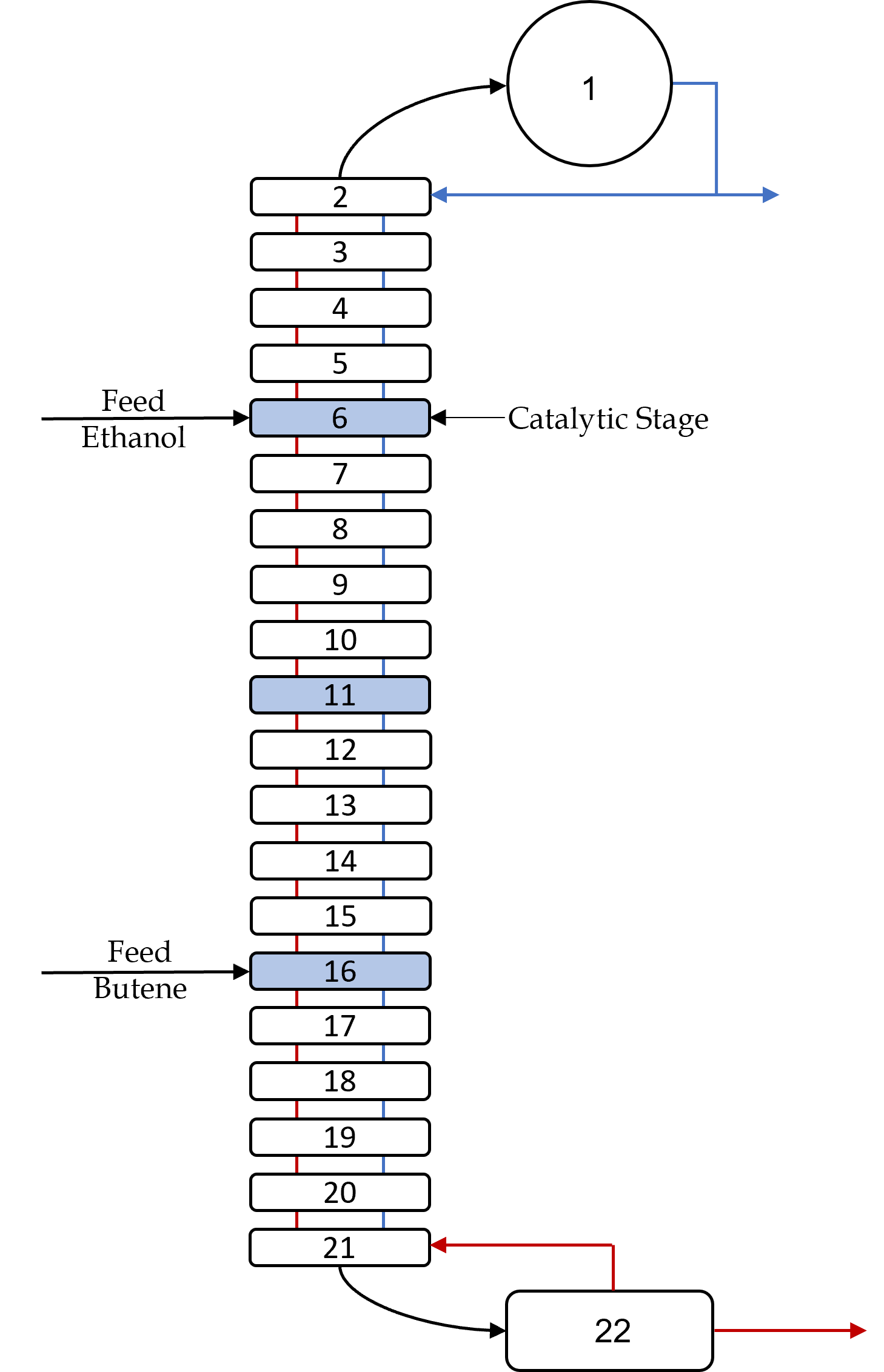}
         \caption{Catalytic column superstructure.}
         \label{fig: catalytic_ratebased}
     \end{subfigure}
     \hfill
     \begin{subfigure}[b]{0.33\textwidth}
         \centering
         \includegraphics[width=\textwidth]{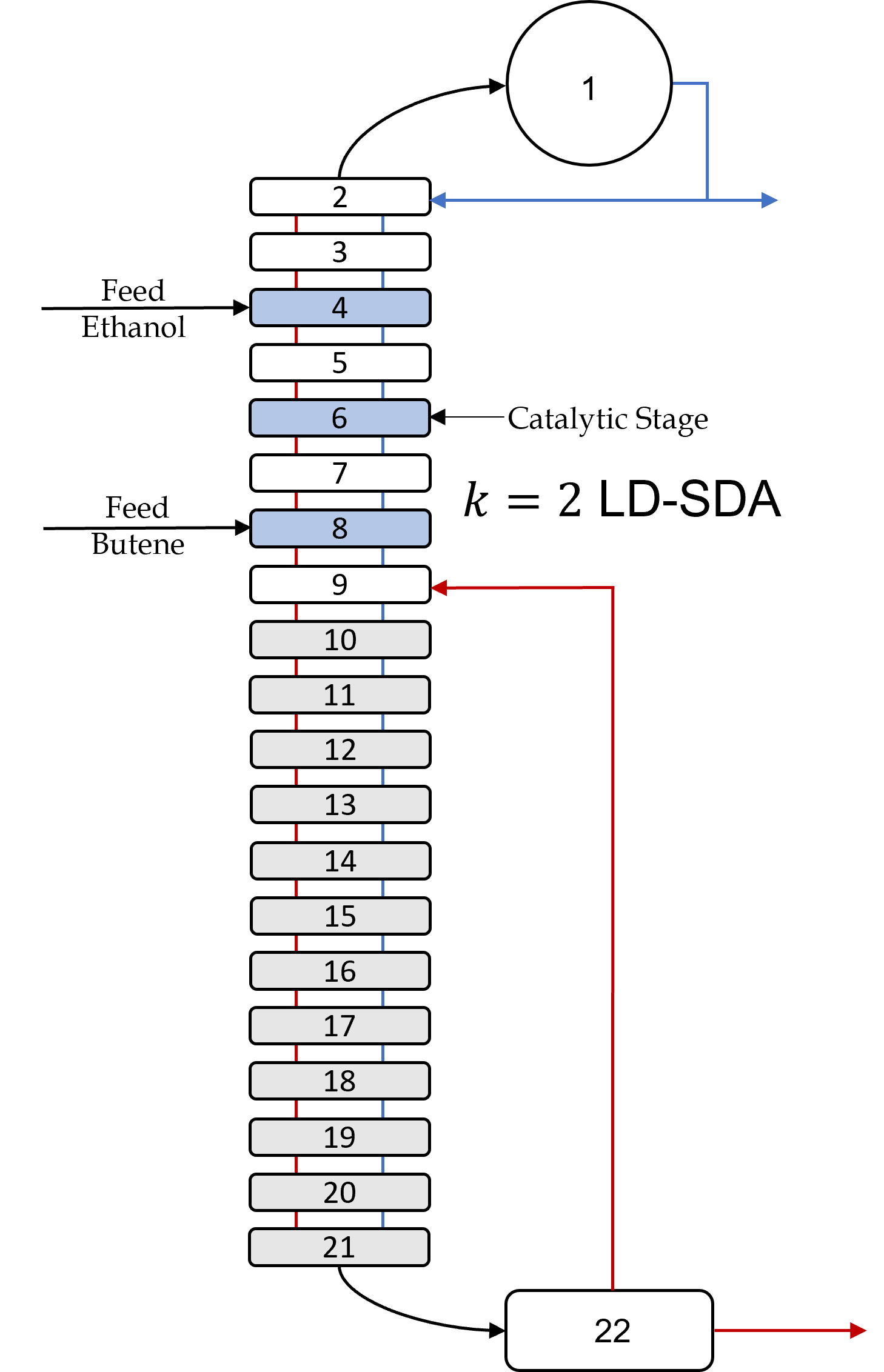}
         \caption{Catalytic distillation column configuration obtained with LD-SDA and $N_2$.}
         \label{fig:catalytic_k2_ldsda}
     \end{subfigure}
     \hfill
     \begin{subfigure}[b]{0.33\textwidth}
         \centering
         \includegraphics[width=\textwidth]{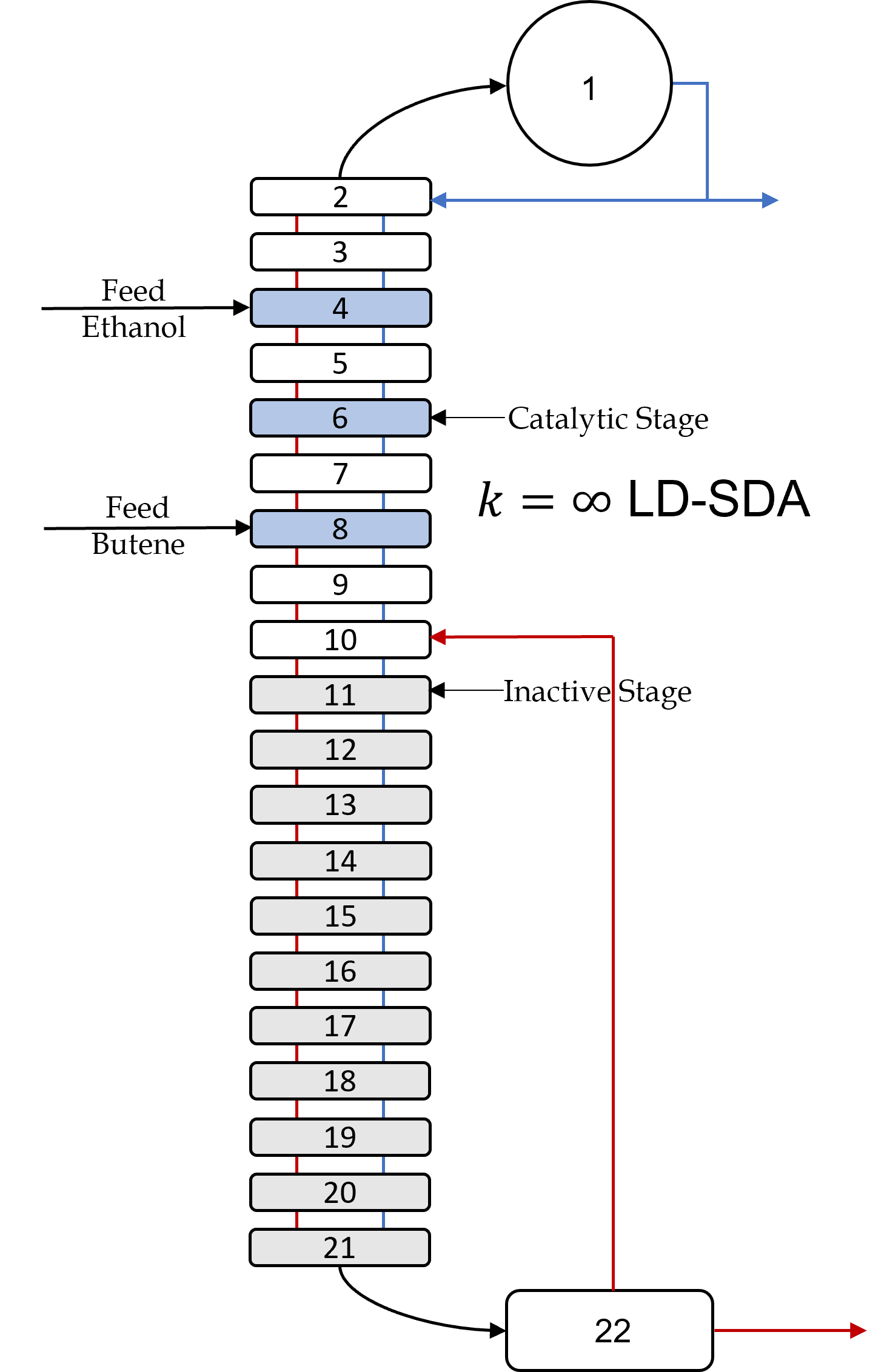}
         \caption{Catalytic distillation column configuration obtained with LD-SDA and $N_\infty$.}
         \label{fig: catalytic_kinf_ldsda}
     \end{subfigure}
        \caption{Visualization of the catalytic distillation columns obtained using LD-SDA with $N_2$ and $N_{\infty}$. The feed trays for Ethanol and Butane are shaded in blue, with existing trays shown in white, while bypass trays are shaded in gray. The boil-up position is marked in red, and the reflux position is depicted in blue.}
        \label{fig:Catalytic Distillation Column }
\end{figure}

\subsection{Optimal Design for Chemical Batch Processing}

Consider an instance of the optimal design for chemical batch processing from~\textcite{kocis1988global} formulated as a GDP.
This is a convexified GDP that aims to find an optimal design for multiproduct batch plants that minimizes the sum of exponential costs.
In our example, the process has three processing stages where fixed amounts of \(q_i\) of two products must be produced.
The goal of the problem is to determine the number of parallel units \(n_j\), the volume \(v_j\) of each stage \(j\), the batch sizes \(b_i\), and the cycle time \(tl_i\) of each product \(i\).
The given parameters of the problem are the time horizon \(h\), cost coefficients \(\alpha_j, \beta_j\) for each stage \(j\), size factors \(s_{ij}\), and processing time \(t_{ij}\) for product \(i\) in stage \(j\).
The optimization model employs Boolean variables \(Y_{kj}\) to indicate the presence of a stage, potentially representing three unit types: mixers, reactors, and centrifuges.
The formulation of the model can be found in Supplementary Material  \ref{app:batch_form}, and the external variable reformulation of the Boolean variables is described in Supplementary Material  \ref{app:batch_reform}.


 
\begin{figure}
    \centering
    \includegraphics[width=0.65\textwidth]{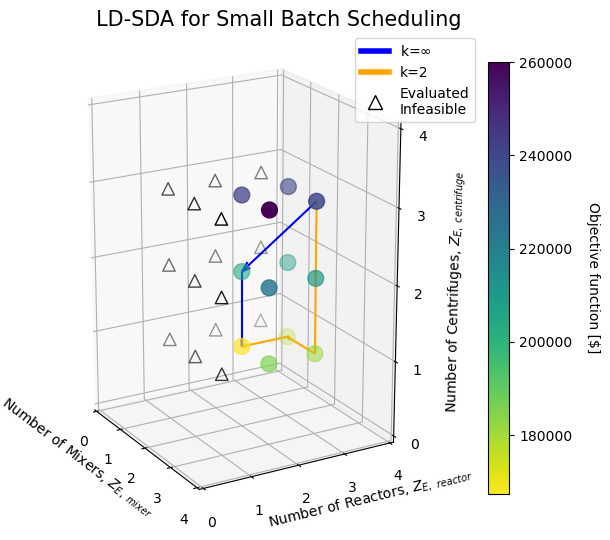}
    \caption{Visualization of the paths traversed by LD-SDA and the solutions found using both neighborhoods for the small batch scheduling case study. The lattice includes points proven to be infeasible with a solver (black triangle). Both solutions were initialized with a configuration containing all mixers, reactors, and centrifuges. LD-SDA with both $k=2$ and $k=\infty$ converged to the same solution $(2,2,1)$, corresponding to a global optimal solution. This example highlights convergence to the same solution in a convex GDP problem using different neighborhoods.}
    \label{fig:small_batch}
\end{figure}

The problem was initialized by setting the maximum number of units, i.e., \((3,3,3)\), for the number of mixers, reactors, and centrifuges, respectively.
The algorithm terminates on a solution with objective \(\$ 167,427\) with external variables \((2,2,1)\) for both $k=2$ or $k=\infty$ neighborhood alternatives in LD-SDA.
The trajectories taken by both searches of the LD-SDA for the small batch problem are shown in Figure~\ref{fig:small_batch}.
This solution corresponds to a global optimal solution of the problem, suggesting that convergence to global optimal solutions in convex GDP might be achieved even with $k=2$ in the Neighborhood Search step.
For this small problem, the solution times were negligible (less than two seconds). 
Still, this example is included to observe convergence to the same solution in a convex GDP problem using different neighborhoods.




\section{Conclusions and Final Remarks}
\label{sec:conclusions}


This work has presented the Logic-Based Discrete-Steepest Descent Algorithm as an optimization method for GDP problems with ordered Boolean variables, which often appear in process superstructure and single-unit design problems. 
The unique characteristics of the LD-SDA are highlighted, and its similarities with other existing logic-based methods are discussed.
To verify the performance of the LD-SDA, we solved various GDP problems with applications in process systems engineering, such as reactor series volume minimization, binary distillation column design, rate-based catalytic distillation column design, and chemical batch process design.
The LD-SDA has demonstrated an efficient convergence toward high-quality solutions that outperformed state-of-the-art MINLP solvers and GDP solution techniques for the problems studied.
The results show that LD-SDA is a valuable tool for solving GDP problems with the special ordered structure considered in this work.
Nonetheless, the scalability of the LD-SDA still needs to be evaluated for larger superstructure problems, e.g., those resulting in more than 7 external variables.
The limitations of the LD-SDA include the lack of guarantee for a globally optimal solution due to its local search nature. 
Additionally, the exponential growth of neighbors with increasing reformulated variables can make neighborhood evaluation prohibitively expensive for large-scale problems.


Future research directions include utilizing the LD-SDA to solve larger and more challenging ordered GDPs.
Similarly, we propose exploring theoretical convergence guarantees of the LD-SDA method, with a special focus on convex GDP problems and their relation to integrally convex problems in discrete analysis.
Moreover, future work also involves the integration of the LD-SDA into the GDPOpt solver in Pyomo.GDP, making this algorithm available to a wider audience.
Finally, we will study the parallelization of NLP solutions in the neighborhood search.
The Neighborhood Search can be faster by dividing the computation involved in solving NLP problems into multiple tasks that can be executed simultaneously, eventually improving the performance of the LD-SDA.

\section*{Acknowledgments}
D.B.N. was supported by the NASA Academic Mission Services, Contract No. NNA16BD14C. D.B.N. and A.L. acknowledge the support of the startup grant of the Davidson School of Chemical Engineering at Purdue University.

\printbibliography
\newpage


\section*{Supplementary Material}
\setcounter{section}{0}
\renewcommand{\thesection}{\Alph{section}} 
\renewcommand{\thesubsection}{\thesection.\arabic{subsection}}

The supplementary material provides detailed formulations of GDP models, including reactors and chemical batch processing systems. 
It outlines the objective functions and presents both the algebraic and logical constraints for each model. 
Additionally, this section explains how Boolean variables are reformulated into external variables. 
Several models are also illustrated through figures for clarity.


\section{Generalized Disjunctive Programming formulations}

This appendix includes the formulations of the examples of the problems solved in this manuscript as Generalized Disjunctive Programs.

\subsection{Series of Reactors}
\label{sec:cstr_formulation}

Set of components (index \textit{i})
\begin{equation}
 I=\{A,B\}
\end{equation}
Set of units in the superstructure (index \textit{n, j})
\begin{equation}
 N=\{1, ... ,NT\} 
\end{equation}

Existence of an unreacted feed in unit \textit{n}
\begin{equation}
\begin{gathered}
 YF_n \in \{\textit{True, False}\} \; \; \forall \ n \in N \\
 \text{If } YF_n=\textit{True} \implies \text{There is unreacted feed in reactor \textit{n}}
 \end{gathered}
\end{equation}

Existence of a recycle flow in unit \textit{n}
\begin{equation}
\begin{gathered}
 YR_n \in \{\textit{True, False}\} \; \;  \forall \ n \in N \\
 \text{If } YR_n=\textit{True} \implies \text{There is recycle in reactor \textit{n}}
 \end{gathered}
\end{equation}

Unit operation in \textit{n}: If at the current unit \textit{n} \textbf{every} unit after it (from one to \textit{n}) \textbf{is not} an unreacted feed \textbf{or} if the current unit \textit{n} has the unreacted feed, then the unit is a CSTR (the opposite is also true)
\begin{equation}
\begin{gathered}
 YP_n \iff \left(\bigwedge_{j \in \{1,2,..n\}}\neg YF_j\right)\vee YF_n \; \;  \forall \ n \in N \\
 \text{If } YP_n=\textit{True} \implies \text{Unit \textit{n} is a CSTR} \\
 \text{If } YP_n=\textit{False} \implies \text{Unit \textit{n} is a bypass}
 \end{gathered}
\end{equation}

The unit must be a CSTR to include a recycle at \textit{n}
\begin{equation}
YR_n \implies YP_n \; \;  \forall \ n \in N
\end{equation}

There is only one unreacted feed
\begin{equation}
\begin{gathered}
    \bigveebar_{n \in N} 
\end{gathered}
YF_n 
\end{equation}

There is only one recycling stream.
\begin{equation}
\begin{gathered}
    \bigveebar_{n \in N} 
\end{gathered}
YR_n 
\end{equation}

Unreacted feed unit: Partial mole balance
\begin{equation}
0=F0_i+ FR_{i,NT}-F_{i,NT}+r_{i,NT} V_{NT} \; \; \forall \ i \in I 
\end{equation}

Unreacted feed unit: Continuity
\begin{equation}
0=Q_{F0}+Q_{FR,NT}-Q_{NT}
\end{equation}

Reactor Sequence: Partial mole balance
\begin{equation}
0=F_{i,n+1}+FR_{i,n}-F_{i,n}+r_{i,n} V_n \; \; \forall \ n \in N \setminus \{NT\},\forall \ i \in I
\end{equation}

Reactor Sequence: Continuity
\begin{equation}
0=Q_{n+1}+Q_{FR,n}-Q_n \; \; \forall \ n \in N \setminus \{NT\}
\end{equation}

If unit \textit{n} is a CSTR or a bypass
\begin{equation}
\left[
\begin{gathered} 
    YP_n \\
    r_{A,n} Q_n^2=-kF_{A,n} F_{B,n} \\
    r_{B,n}=-r_{A,n}\\
    c_n = V_n
    \end{gathered}
    \right]
    \bigveebar
    \left[
\begin{gathered} 
    \neg YP_n \\
    FR_{i,n}=0  \; \; \forall \ i \in I\\
    r_{i,n}=0  \; \; \forall \ i \in I\\
    Q_{FR,n}=0  \\
    c_n=0
    \end{gathered}
    \right]
    \forall \ n \in N
\end{equation}

If there is recycle in before reactor \textit{n}
\begin{equation}
\left[
\begin{gathered} 
    YR_n \\
    FR_{i,n}=R_i \; \; \forall \ i \in I\\
    Q_{FR,n}=Q_R\\
    \end{gathered}
    \right]
    \bigveebar
    \left[
\begin{gathered} 
    \neg YR_n \\
    FR_{i,n}=0  \; \; \forall \ i \in I\\
    Q_{FR,n}=0  \\
    \end{gathered}
    \right]
    \forall \ n \in N
\end{equation}

Splitting point: Partial mole balance
\begin{equation}
0=F_{i,1}-P_i-R_i \; \; \forall \ i \in I 
\end{equation}

Splitting point: Continuity
\begin{equation}
0=Q_1-Q_P-Q_R 
\end{equation}

Splitting point: Additional constraint
\begin{equation}
0=P_i Q_1-F_{i,1} Q_P \; \; \forall \ i \in I 
\end{equation}

Product specification constraint
\begin{equation}
0.95Q_P=P_B 
\end{equation}

Volume constraint 
\begin{equation}
V_n=V_{n-1} \; \; \forall \ n \in N \setminus \{1\}
\end{equation}

Objective Function: Total reactor network volume
\begin{equation}
 f_{OBJ} = \min \sum_{n \in N}{c_n}
\end{equation}

\subsection{Distillation Column Design}
\label{sec:column_formulation}
Set of trays (index \textit{t})
\begin{equation}
T = \{2,3, \dots, 16\}
\end{equation}

Set of composition (index \textit{c})
\begin{equation}
C = \{\text{Benzene, Toluene} \}
\end{equation}

Existence of tray \textit{t}
\begin{equation}
\begin{gathered}
Y_t \in \{\textit{True, False}\} \; \;  \forall \ t \in T \\
\text{If } Y_i=\textit{True} \implies \text{There exist a tray in stage \textit{t}}
\end{gathered}
\end{equation}

Existence of boil-up flow in tray \textit{t}
\begin{equation}
\begin{gathered}
YB_t \in \{\textit{True, False}\} \; \;  \forall \ t \in T \\
\text{If } YB_t=\textit{True} \implies \text{There is a boil-up flow in tray \textit{t}}
\end{gathered}
\end{equation}

Existence of reflux flow in tray \textit{t}
\begin{equation}
\begin{gathered}
YR_t \in \{\textit{True, False}\} \; \;  \forall \ t \in T \\
\text{If } YR_i=\textit{True} \implies \text{There is a reflux flow in tray \textit{t}}
\end{gathered}
\end{equation}

There is only one boil-up flow in the distillation column.
\begin{equation}
\begin{gathered}
    \bigveebar_{t \in T} 
\end{gathered}
YB_t 
\end{equation}

There is only one reflux flow in the distillation column.
\begin{equation}
\begin{gathered}
    \bigveebar_{t \in T}
\end{gathered}
YR_t 
\end{equation}

Tray \(t\) is an equilibrium stage or a bypass, where \(\mathbf{g}_1 (\mathbf{x})\) contains equilibrium mass and energy balances, while \(\mathbf{g}_2 (\mathbf{x})\) contains a bypass material balance
\begin{equation}
\begin{bmatrix}
Y_t \\
\mathbf{g}_1 (\mathbf{x}) = \mathbf{0} \\ 
y_{t,\text{active}}=1 
\end{bmatrix}
\bigveebar
\begin{bmatrix}
\neg Y_t \\
\mathbf{g}_2 (\mathbf{x}) = \mathbf{0} \\ 
y_{t,\text{active}}=0 
\end{bmatrix},
\quad \forall \ t \in T
\end{equation}

If the reflux flow is on or above tray \(t\) and the boil-up flow is on or below tray \(t\), then tray \(t\) is an equilibrium stage (the opposite is also \(True\)).

\begin{equation}
\left(
\begin{gathered}
    {\bigvee} \\
    \forall \ \tau \in \{ t, t+1,\dots,16 \} 
\end{gathered}
YR_{\tau}
\right)
\land
\left(
\begin{gathered}
    {\bigwedge} \\
    \forall \ \tau \in \{ t, t+1,\dots,16 \} 
\end{gathered}
\neg YB_{\tau} \vee YB_t
\right)
\Longleftrightarrow
Y_t,
\quad \forall \ t \in T
\end{equation}

The reflux flow stage is not below the feed tray.
\begin{equation}
\begin{gathered}
    {\bigwedge} \\
    \forall \ t \in \{ 2,3,\dots,8 \} 
\end{gathered}
\neg YR_t
\end{equation}

The boil-up flow is not above the feed tray.
\begin{equation}
\begin{gathered}
    {\bigwedge} \\
    \forall \ t \in \{ 10,11,\dots,17 \} 
\end{gathered}
\neg YB_t
\end{equation}

The column has at least eight active trays.
\begin{equation}
\sum_{t \in T} y_{t,active} \geq 8
\end{equation}

Tray 1 (reboiler), tray 9 (feed), and tray 17 (condenser) are equilibrium stages.
\begin{equation}
Y_1 = True, \ Y_9 = True, \ Y_{17} = True,
\end{equation}

Benzene concentration constraint (distillate product)
\begin{equation}
X_{D, \text{Benzene}} \geq 0.95
\end{equation}

Toluene concentration constraint (bottom product)
\begin{equation}
X_{B, \text{Toluene}} \geq 0.95
\end{equation}

Bounds imposed over the reflux ratio \(RF\)
\begin{equation}
0.5 \leq RF \leq 4
\end{equation}

Bounds imposed over the Reboil ratio \(RB\)
\begin{equation}
1.3 \leq RB \leq 4
\end{equation}

Objective Function: Sum of the capital cost (number of active trays) and operating cost (reboiler and condenser duty), where \(Q_R\) is the reboiler duty, \(Q_C\) is the condenser duty, and the sum term over binary variables \(y_{t, \text{active}} \) represents the total number of active trays.

\begin{equation}
f_{OBJ} = \min_{RF, RB, feedtray, N_{trays}} 10^3 \left( (Q_R + Q_C) + \sum_{t \in T}^{N_{trays}} y_{t, \text{active}} \right)
\end{equation}

\subsection{Small Batch Problem}
\label{app:batch_form}
Set of components (index \textit{i})
\begin{equation}
 I=\{A,B\}
\end{equation}

Set of stages (index \textit{j})
\begin{equation}
 J=\{\text{mixer}, \text{reactor}, \text{centrifuge}\}
\end{equation}

Set of potential number of parallel units for each stage (index \textit{k})
\begin{equation}
    K = \{1, 2, 3\}
\end{equation}

Existence of the parallel units for each stage \textit{j}
\begin{equation}
\begin{gathered}
 Y_{kj} \in \{\textit{True, False}\} \; \forall \ k \in K, \ j \in J \\
 \text{If } Y_{kj} =\textit{True} \implies \text{There are \textit{k} parallel units in stage \textit{j}}
 \end{gathered}
\end{equation}

Only one of the parallel unit existence is \(True\)
\begin{equation}
\begin{gathered}
    \bigveebar_{j \in J} 
\end{gathered}
Y_{kj}  \quad \forall \ k \in K
\end{equation}

Volume requirement in stage \textit{j}
\begin{equation}
    v_j \geq \ln (s_{ij}) + b_i \quad \forall \ i \in I, \ j \in J 
\end{equation}

Cycle time for each product \textit{i}
\begin{equation}
   n_j + t_{Li} \geq \ln (t_{ij}) \quad \forall \ i \in I, \ j \in J 
\end{equation}

Constraint for production time(horizon constraint)
\begin{equation}
    \sum_{i \in I} Q_i \exp{(t_{Li} - b_i)} \leq H
\end{equation}

Relating number of units to \(0-1\) variables
\begin{equation}
    n_j = \sum_{k \in K} \gamma_{kj} \quad \forall \ j \in J 
\end{equation}

If only \textit{k} parallel units exist in stage \textit{j}

\begin{equation}
\begin{bmatrix}
    \begin{gathered}
    Y_{kj} \\
    \gamma_{kj} = \ln{(k)} \\
    \end{gathered}
\end{bmatrix}
\bigveebar
\begin{bmatrix}
    \begin{gathered}
    \neg Y_{kj} \\
    \gamma_{kj} = {0}\\
    \end{gathered}
\end{bmatrix}
\quad \forall \ k \in K, \ j \in J
\end{equation}

Objective Function: the investment cost for setting the small batch system \([\$]\)

\begin{equation}
f_{OBJ} = \min \sum_{j \in J} \alpha_j ( \exp(n_j + \beta_j v_j))
\end{equation}

\section{External variable reformulation for example problems}

This appendix presents the external variable reformulation of the Boolean variables in the examples considered in this manuscript.

\subsection{Series of Reactors Problem}
\label{app:cstr_reform}
\numberwithin{equation}{subsection}

\begin{equation}
YF_n = 
    \begin{cases}
    True, \;\; \mathbf{z}_{E,1}=n \\
    False, \;\; \text{otherwise} 
    \end{cases}
    \forall \ n \in N
\end{equation}

\begin{equation}
YR_n = 
    \begin{cases}
    True, \;\; \mathbf{z}_{E,2}=n \\
    False, \;\; \text{otherwise} 
    \end{cases}
    \forall \ n \in N
\end{equation}

\begin{equation}
X_1 = 
\left\{
    \mathbf{z_E} \in \mathbb{Z} ^2:
    \begin{gathered}
        1 \leq \mathbf{{z}_{E,1}} \leq NT \\
        1 \leq \mathbf{{z}_{E,2}} \leq NT
    \end{gathered}
\right\}
\end{equation}

\begin{equation}
X_2 = 
\left\{
    \mathbf{{z}_E} \in \mathbb{Z} ^2: \mathbf{{z}_{E,2}}-\mathbf{z_{E,1}} \leq 0
\right\}
\end{equation}

\begin{equation}
X = X_1 \cap X_2 =
\left\{
    x_E \in \mathbb{Z} ^2:
    \begin{gathered}
        1 - \mathbf{{z}_{E,1}} \leq 0 \\
        \mathbf{{z}_{E,1}} - NT \leq 0 \\
        1 - \mathbf{{z}_{E,2}} \leq 0 \\
        \mathbf{{z}_{E,2}} - \mathbf{{z}_{E,1}} \leq 0 \\
    \end{gathered}
\right\}
\end{equation}

\subsection{External Variable Reformulation for Distillation Column Problem}
\label{app:column_reform}

\begin{figure}[H]
    \centering
    \includegraphics[width=0.75\textwidth]{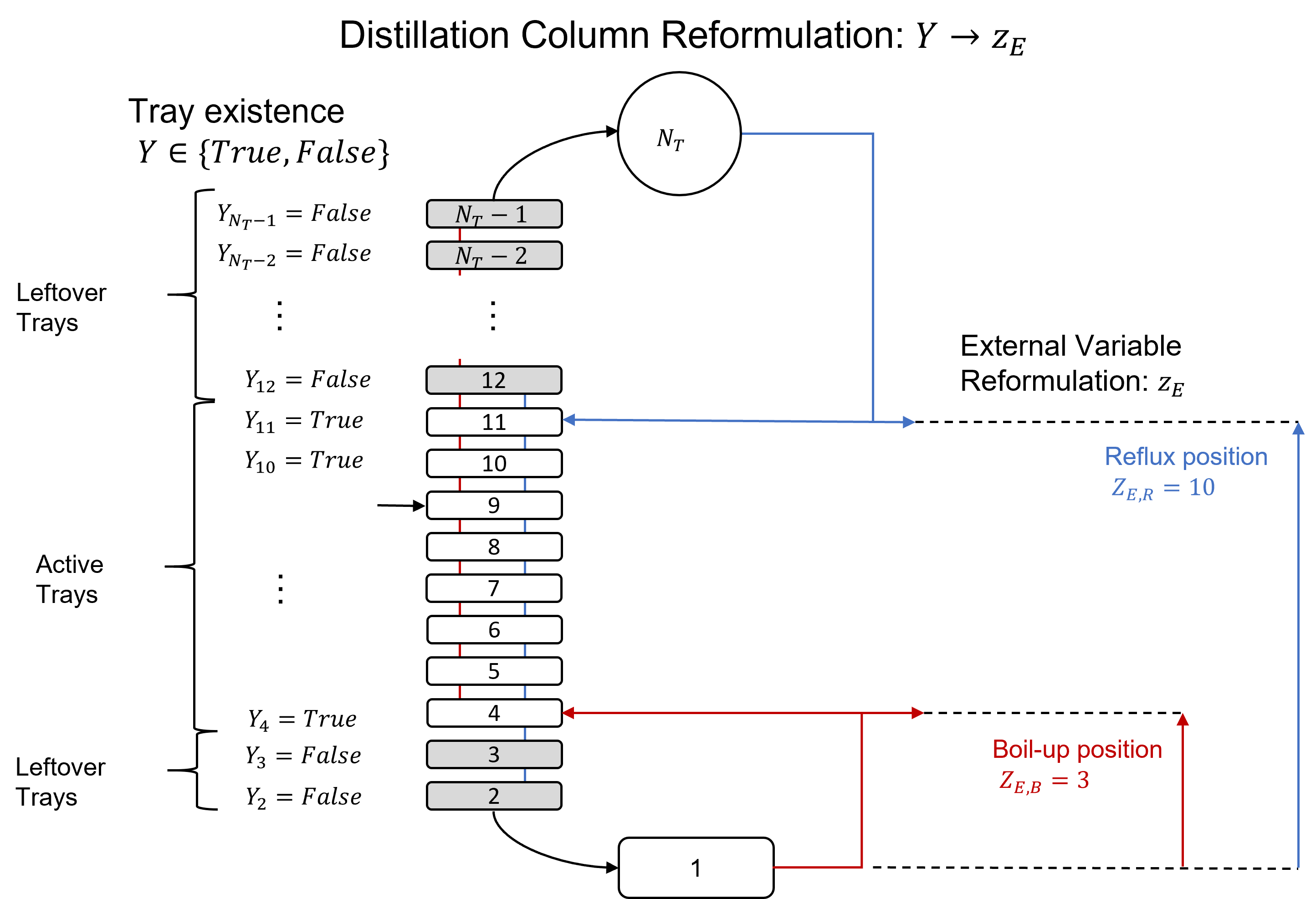}
    \caption{Example of distillation column external variable reformulation}
    \label{fig:column_reformulation}
\end{figure}

\begin{equation}
YR_t = 
\begin{cases}
    True, \;\; \mathbf{z}_{\mathbf{E},\text{reflux}}=t-1 \\
    False, \;\; \text{otherwise} 
    \end{cases}
    \forall \ t \in T
\end{equation}

\begin{equation}
YB_t = 
\begin{cases}
    True, \;\; \mathbf{z}_{\mathbf{E},\text{boil-up}}=t-1 \\
    False, \;\; \text{otherwise} 
    \end{cases}
    \forall \ t \in T
\end{equation}

\begin{equation}
X_1 = 
\left\{
    \mathbf{z}_\mathbf{E} \in \mathbb{Z} ^2:
    \begin{gathered}
        1 \leq \mathbf{z}_{\mathbf{E},\text{reflux}} \leq 15 \\
        1 \leq \mathbf{z}_{\mathbf{E},\text{boil-up}} \leq 15 \\
    \end{gathered}
\right\}
\end{equation}

\subsection{Catalytic Distillation Column Problem}

\begin{figure}[H]
    \centering
    \includegraphics[width=0.60\textwidth]{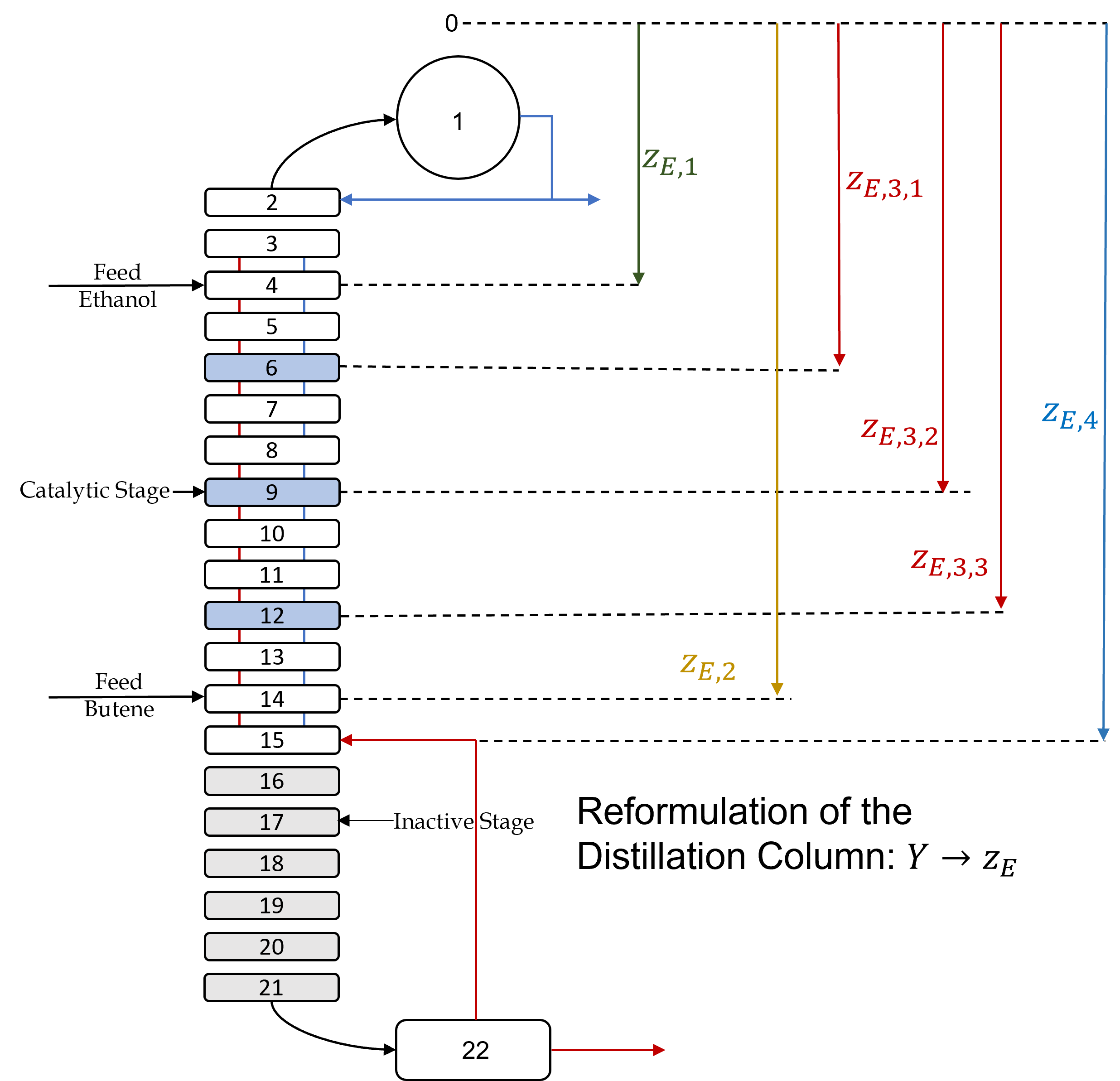}
    \caption{Example of catalytic distillation column external variable reformulation}
    \label{fig:catalytic_reformulation}
\end{figure}

The external variable reformulation is equivalent to the one presented in \parencite{linan2020optimalII} with Boolean variables instead of binary variables.
We highlight in Figure~\ref{fig:catalytic_reformulation} how these external variables are interpretable as relative positions of the ethanol feed, the butene feed, the catalytic stages and the boil-up.

\subsection{Small Batch Problem}
\label{app:batch_reform}
\begin{equation}
Y_{k,\text{mixer}} = 
    \begin{cases}
    True, \;\; \mathbf{z}_{\mathbf{E},\text{mixer}}=k \\
    False, \;\; \text{otherwise} 
    \end{cases}
    \forall \ k \in K
\end{equation}

\begin{equation}
Y_{k,\text{reactor}} = 
    \begin{cases}
    True, \;\; \mathbf{z}_{\mathbf{E},\text{reactor}}=k \\
    False, \;\; \text{otherwise} 
    \end{cases}
    \forall \ k \in K
\end{equation}

\begin{equation}
Y_{k,\text{centrifuge}} = 
    \begin{cases}
    True, \;\; \mathbf{z}_{\mathbf{E},\text{centrifuge}}=k \\
    False, \;\; \text{otherwise} 
    \end{cases}
    \forall \ k \in K
\end{equation}

\begin{equation}
X_1 = 
\left\{
    \mathbf{z}_\mathbf{E} \in \mathbb{Z} ^3:
    \begin{gathered}
        1 \leq \mathbf{z}_{\mathbf{E},\text{mixer}} \leq K \\
        1 \leq \mathbf{z}_{\mathbf{E},\text{reactor}} \leq K \\
        1 \leq \mathbf{z}_{\mathbf{E},\text{centrifuge}} \leq K
    \end{gathered}
\right\}
\end{equation}


\end{document}